\documentclass[3p, 11pt]{article}%

\usepackage{titling}
\usepackage[hmarginratio=1:1,top=25mm,left=25mm,columnsep=1cm]{geometry} %
\usepackage{array}
\usepackage{color}
\usepackage{tabularx}
\usepackage{graphicx}
\usepackage[svgnames,table,xcdraw]{xcolor}
\usepackage{amsmath}
\usepackage{amsthm}
\usepackage{amssymb}
\usepackage{comment}
\usepackage{appendix}
\usepackage{dsfont}
\usepackage{bm}
\usepackage{multirow}
\usepackage{fourier} 
\usepackage{epstopdf}
\usepackage[font=small,labelfont=bf]{caption}	
\usepackage{multicol,booktabs}
\usepackage{pbox} 
\usepackage{subcaption} 

\graphicspath{{../}}	

\usepackage{hyperref}
\hypersetup{
    colorlinks=true,
    linkcolor=Blue, 
    citecolor=DarkRed, 
    urlcolor=blue  } 

\newcommand{\x}{\mbf{x}}

\newcommand{\Th}{\mathcal{T}_h}				

\newcommand{\mc}[1]{\mathcal{#1}}			
\newcommand{\mbf}[1]{\mathbf{#1}}			%

\newcommand{\m}{\mathbf{m}}

\newcommand{\f}{\mathbf{f}}

\newcommand{\ddd}[1]{\text{d} #1}

\newcommand{\bdm}{\begin{displaymath}}
    \newcommand{\edm}{\end{displaymath}}

\newcommand{\bea}{\begin{eqnarray} }
    \newcommand{\eea}{\end{eqnarray} }

\renewcommand{\div}{{\nabla \cdot}}
\newcommand{\divh}{{\nabla_h \cdot}}

\newcommand{\HH}{{\mathbf{H}}}

\newcommand{\rr}{{\mathbf{r}}}
\newcommand{\ee}{{\mathbf{e}}}

\newcommand{\vv}{{\mathbf{v}}}

\newcommand{\mm}{{\mathbf{m}}}

\newcommand{\PP}{{\mathbf{P}}}

\newcommand{\pp}{{\mathbf{p}}}
\newcommand{\qq}{{\mathbf{q}}}
\newcommand{\nn}{{\mathbf{n}}}

\newcommand{\dhcentral}{\partial_{h}^{s,c}}
\newcommand{\dhp}{\partial_{h}^{s,+}}
\newcommand{\dhm}{\partial_{h}^{s,-}}

\newcommand{\jump}[1]{\bigl\llbracket{#1}\bigr\rrbracket}
\newcommand{\mean}[1]{\overline{\left({#1}\right)}}
\newcommand{\phibf}{\bm{\varphi}}
\newcommand{\pK}{\partial K}
\newcommand{\DMV}{\mathcal{V}_{t,x}}

\newcommand{\tvr}{\tilde \rho}

\newcommand{\tvm}{{\tilde \m}}

\newcommand{\tvS}{\tilde S}
\newcommand{\Td}{\Omega}
\newcommand{\RR}{\Bbb R}
\newcommand{\Div}{{\nabla \cdot}}
\newcommand{\Grad}{\nabla}
\newcommand{\bfphi}{\boldsymbol{\varphi}}

\newcommand{\dx}{\,{\rm d} {x}}
\newcommand{\dt}{\,{\rm d} t }

\newcommand{\intTd}[1]{\int_{\Omega} #1 \ \dx}

\newtheorem{theorem}{Theorem}[section]
\newtheorem{lemma}[theorem]{Lemma}

\newtheorem{remark}[theorem]{Remark}
\newtheorem{assumption}[theorem]{Assumption}
\newtheorem{definition}[theorem]{Definition}

\usepackage{letltxmacro}

\LetLtxMacro\davidsincludegraphics\includegraphics

\makeatletter

\def\@starttocbutdonotshowit#1{%
    \begingroup
    \makeatletter
    \if@filesw
    \expandafter\newwrite\csname tf@#1\endcsname
    \immediate\openout \csname tf@#1\endcsname \jobname.#1\relax
    \fi
    \@nobreakfalse
    \endgroup}

\newcommand{\listoffigurenumbernames}{%
    \@starttocbutdonotshowit{lfn}%
}

\date{}

\begin{document}

    \title{
        \textbf{Convergence of a Hyperbolic Thermodynamically Compatible Finite Volume scheme for the Euler equations}
    }

    \author{
        Michael Dumbser\thanks{Department of Civil, Environmental and Mechanical Engineering, University of Trento, Via Mesiano 77, I-38123 Trento, Italy,  (\href{mailto:michael.Dumbser@unitn.it}{michael.dumbser@unitn.it})
        },\
        \quad
        M{\'a}ria Luk{\'a}{\v{c}}ov{\'a}-Medvid'ov{\'a},\thanksgap{-0.5ex}\thanks{Institut f\"ur Mathematik, Johannes-Gutenberg-Universit\"at Mainz, Staudingerweg 9, 55099 Mainz, Germany,
            (\href{mailto:lukacova@mathematik.uni-mainz.de}{lukacova@mathematik.uni-mainz.de})}
        \quad
        and
        \quad
        Andrea Thomann\thanksgap{-0.5ex}
        \thanks{Universit\'e de Strasbourg, CNRS, Inria, IRMA, F-67000 Strasbourg, France,
            (\href{mailto:andrea.thomann@inria.fr}{andrea.thomann@inria.fr})}
        \textsuperscript{~,$\ast$}
    }
    \thanksmarkseries{arabic}

    \maketitle 

    {\let\thefootnote\relax\footnote{\hspace{0.25cm}$^* $Corresponding author}}

    \paragraph{Abstract}

    We study the convergence of a novel family of thermodynamically compatible schemes for hyperbolic systems (HTC schemes) in the framework of dissipative weak solutions, applied to the Euler equations of compressible gas dynamics. Two key novelties of our method are i) entropy is treated as one of the main field quantities and ii) the total energy conservation is a consequence of compatible discretization and application of the Abgrall flux.




\section{Introduction} \label{sec:introduction}

In his seminal work \cite{Godunov1961} Godunov found the connection between symmetric hyperbolicity  in the sense of Friedrichs \cite{FriedrichsSymm} and thermodynamic compatibility. A similar connection was later also established by Friedrichs \& Lax in \cite{FriedrichsLax}. A key ingredient of the work of Godunov was the assumption of an underlying variational structure of the governing PDE system, which is reasonable for many physical systems. Indeed, Godunov showed that for hyperbolic systems which have an underlying variational formulation, the total energy conservation equation is an \textit{extra conservation law} that can be obtained as a consequence of the other equations by taking the dot product of the original system with the so-called thermodynamic dual variables, which are the partial derivatives of the total energy potential with respect to the conservative variables.

Important contributions to the topic of symmetric hyperbolic and thermodynamically compatible (SHTC) systems were subsequently made by Boillat \cite{Boillat74} and Ruggeri \cite{Ruggeri81}. In the latter reference the thermodynamic dual variables were also called the \textit{main field}, while other works refer to them as the \textit{Godunov variables}, see e.g. \cite{Freistuehler2019}.
In later work by Godunov \& Romenski and collaborators, the original theory of Godunov was extended to a much larger class of hyperbolic systems, including magnetohydrodynamics (MHD), nonlinear hyperelasticity, compressible multi-phase flows and even relativistic gasdynamics, see e.g.
\cite{God1972MHD,GodunovRomenski72,Rom1998,RomenskiTwoPhase2010,Godunov2012}.
A connection between the SHTC formalism and Hamiltonian mechanics was recently found in \cite{SHTC-GENERIC-CMAT}.

Usually in the context of SHTC systems, the \textit{entropy} is treated as \textit{main evolution} quantity, while the \textit{total energy conservation law} is the extra conservation law that is obtained as a \textit{consequence}. The privileged role of the total energy is motivated by the underlying variational formulation.

Up to now, most thermodynamically compatible schemes fall into the well-known framework established by Tadmor \cite{Tadmor1}, in which discrete compatibility with the \textit{entropy inequality} is obtained as a consequence of a compatible discretization of all other equations. Without pretending completeness, an overview of important work made in the field can be found in \cite{ FjordholmMishraTadmor,CastroFjordholm,GassnerEntropyEuler,ShuEntropyMHD1,ChandrashekarKlingenberg2016,Ranocha2020}.

For scalar conservation laws the convergence analysis of entropy stable schemes via compensated compactness arguments was shown in \cite{Chatterjee2020}. Although the weak entropy solutions are unique for the scalar hyperbolic conservation laws, it is not true for multidimensional hyperbolic systems. In particular, it is well-known that multidimensional Euler equations of gas dynamics are ill-posed in the class of weak entropy solutions. Consequently, a new strategy based on the generalized, the so-called dissipative weak solutions, has been proposed in order to study the convergence of entropy stable numerical methods in \cite{FeiLukMiz2019}, see also \cite{Feireisl2021book, FLM2020a, ALO2023, LO2023}.  There is a \emph{striking similarity between dissipative weak solutions and SHTC systems}: the main field equations are the conservation of mass, momentum and entropy inequality, while the energy balance is a consequence of the latter. We note that in general, only the total energy dissipation holds for dissipative weak solutions. However, if a strong solution to the Euler system exists, the total energy conservation holds.

A general and at the same time very simple framework for the construction of compatible numerical methods that satisfy extra energy conservation laws exactly at the discrete level was recently introduced by Abgrall and collaborators, see e.g.  \cite{Abgrall2018, AbgrallOeffnerRanocha, Abgrall2023HTC}. This framework will also be a basis of the numerical scheme studied in this paper.

An attempt to achieve discrete total energy conservation as a consequence of a compatible discretization of all other equations was made in \cite{SWETurbulence,Busto2022,HTCMHD,Abgrall2023HTC,Thomann2023,HTCLagrange}, but a rigorous convergence analysis of these schemes has not yet been carried out.
The main goal of the present paper is to provide a rigorous convergence analysis of this new class of schemes for the compressible Euler equations via dissipative weak solutions under the assumption of a uniform lower bound on the density and upper bound on the energy.

The rest of this paper is organized as follows: in Sections \ref{sec.model} and \ref{sec:NumericalScheme} we present the governing equations, their vanishing viscosity regularization and the numerical method to be studied in this work. In Sections \ref{sec:stability}, \ref{sec:Consistency} and \ref{sec:convergence} we analyze the stability, consistency and convergence of the scheme. The paper closes with some concluding remarks and an outlook to future research in Section \ref{sec:conclusions}.

    \setcounter{footnote}{0}
    \section{The compressible Euler equations}
    \label{sec.model}

    \subsection{Inviscid equations}

    In this work, we consider the complete Euler equations in the following formulation
    \begin{subequations}
        \label{eq:Euler_in_S}
    \begin{eqnarray}
        && \partial_t \rho + \div (\rho \vv) = 0, \\
        && \partial_t (\rho \vv) + \div (\rho \vv \otimes \vv) + \nabla p(\rho,s) = 0, \\
        && \partial_t (\rho s) + \div (\rho s \vv) \geq 0,
    \end{eqnarray}
    \end{subequations}
    where $\rho = \rho(t,x)$ denotes the mass density, $\vv = \vv(t,x)$ the velocity field and $s = s(t,x)$ the specific entropy.
    In addition, we have the following conservation law for the total energy $\rho E$ given by
    \begin{equation}
        \label{eq:Energy_conservation}
        \partial_t (\rho E) + \div ((\rho E + p) \vv) = 0,
    \end{equation}
    where the specific total energy $E$ is given by the contributions of the internal and kinetic energies
    \begin{equation}
        E = e + \frac{1}{2} {|}\vv {|}^2.
    \end{equation}
    Therein, $e(\rho,s)$ denotes the specific internal energy, which is given by an ideal gas law as follows
    \begin{equation}
        \label{eq:specific_energy_EOS}
        e(\rho,s) = \frac{\rho^{\gamma-1}}{\gamma -1} \exp\left(\frac{s}{c_v}\right).
    \end{equation}
    Therein, $\gamma > 1$ denotes the ratio of specific heats and $c_v$ the specific heat at constant volume.
    Using \eqref{eq:specific_energy_EOS}, the pressure $p$ and temperature $T$ are defined as follows
    \begin{equation}
        \label{eq:pressure_EOS}
       p(\rho,s) =\rho^2\frac{\partial e}{\partial \rho} = (\gamma - 1) \rho e(\rho,s), \quad T = \frac{\partial e}{\partial s} = \frac{e(\rho,s)}{c_v}.
    \end{equation}
    System \eqref{eq:Euler_in_S} is hyperbolic and exhibits the following waves in normal flow direction $\vv_n$
    \begin{equation}
        \label{eq:Eigenvalues}
        \lambda_1 = \vv_n - c, \quad \lambda_{2,\dots,d+1} = \vv_n, \quad \lambda_{d+2} = \vv_n+c.
    \end{equation}
    Therein, $c$ denotes the sound speed given by
    \begin{equation}
        c = \sqrt{\gamma \frac{p}{\rho}}.
    \end{equation}
    Later on, we will use the so called dual variables, which are defined by $\pp = \nabla_\qq (\rho E) \in \mathbb{R}^{d+2}$, where $\qq \in \mathbb{R}^{d+2}$ denotes the vector of state variables $\qq = (\rho, \rho \vv, \rho s)^T \equiv (\rho, \m, S)^T$.  Here, we have used the notation $\mm = \rho \vv$ for the momentum and $S = \rho s$  for the total entropy.
    The respective entries in $\pp$ are given by
    \begin{subequations}
        \label{eq:QtoP}
    \begin{align}
        p_1(\qq) &= \frac{\partial (\rho E)}{\partial \rho} = \left(\gamma-\frac{S}{\rho c_v}\right)\frac{\rho^{\gamma-1}}{\gamma - 1} \exp\left(\frac{S}{\rho c_v}\right) - \frac{1}{2\rho^2}{|}\m {|}^2 =: r,\\
        p_{k+1}(\qq) &= \frac{\partial (\rho E)}{\partial (\rho v_k)}= \frac{\rho v_k}{\rho} = v_k, \\
        p_{d+2}(\qq) &=\frac{\partial (\rho E)}{\partial S} = \frac{1}{c_v} \frac{\rho^{\gamma-1}}{\gamma -1} \exp\left(\frac{S}{\rho c_v}\right) = T,
    \end{align}
    \end{subequations}
    where $k = 1, \dots, d$ denotes the vector index.
    Thus, summarizing, we have $\pp = (r, \vv, T)^T$.
    The transform \eqref{eq:QtoP} from $\qq$ to $\pp$ will be denoted by $\mathfrak{T}^{-1}: \qq \mapsto \pp$.
    Vice versa, the transform from $\pp$ to $\qq$ denoted by $\mathfrak{T}: \pp \mapsto \qq$ is given by
    \begin{subequations}
    \label{eq:PtoQ}
    \begin{align}
        q_1(\pp) &= \exp \left(\frac{- \gamma  c_v T+ c_v T \log \left((\gamma -1)c_v T\right)+r+{|}\vv{|}^2/2}{(\gamma -1) c_v T}\right) = \rho ,\\
        q_{k+1}(\pp) &= v_k \rho = m_k, \\
        q_{d+2}(\pp) &= \rho c_v \log\left(c_v (\gamma - 1)\frac{T}{\rho^{\gamma - 1}}\right) = S.
    \end{align}
    \end{subequations}
    Further, we denote the Jacobian of $\mathfrak{T}^{-1}(\qq)$ by $ \nabla_\qq \mathfrak{T}^{-1}(\qq) =: \HH = \nabla_\qq^2 (\rho E)$, which is the Hessian of the total energy $\rho E$.
    In \cite[Lemma 7.3]{Feireisl2021book} it was shown, that the total energy $\rho E$  as a function of $(\rho,\mm,S)$ employing an ideal gas law is strictly convex for $\rho > 0$.
    Thus $\HH$ is positive definite and invertible.


    \subsection{Vanishing viscosity regularization}
    Let us rewrite the above given set of equations \eqref{eq:Euler_in_S} in a vanishing viscosity regularisation with $\varepsilon > 0$ in the following form
    \begin{subequations}
        \label{eq:Euler_in_S_vvr}
        \begin{eqnarray}
            && \partial_t \rho + \div (\rho \vv) = \div \left(\varepsilon \nabla \rho\right), \label{eq:Euler_in_S_vvr:dens} \\
            && \partial_t (\rho \vv) + \div (\rho \vv \otimes \vv) + \nabla p(\rho,s) = \div\left(\varepsilon \nabla \rho\vv\right), \\
            && \partial_t S + \div (S \vv) = \div\left(\varepsilon\nabla S\right) + \Pi(\qq,\nabla\qq) \geq 0. \label{eq:Euler_in_S_vvr:entropy}
        \end{eqnarray}
    \end{subequations}
    Note that therein, we have rewritten the entropy inequality into a balance law introducing a source term $\Pi \geq 0$ which will be determined in consistency with the energy conservation.
    Analogously, we write the regularized equation for the total energy
    \begin{equation}
        \label{eq:Energy_conservation_vvr}
        \partial_t (\rho E) + \div ((\rho E + p) \vv) = \div\left(\varepsilon\nabla \rho E\right).
    \end{equation}
    Then, \eqref{eq:Euler_in_S_vvr} can be written in the general compact form
    \begin{subequations}
        \label{eq:Euler_compact}
        \begin{eqnarray}
            \partial_t \qq + \div \f(\qq) = \div\left(\varepsilon\nabla\qq\right)+ \mathbf{P(\qq,\nabla \qq)} && \text{in } \Omega \times (0,T)\\
            \qq(0,\cdot~) = \qq^0, && \text{in } \Omega
        \end{eqnarray}
    \end{subequations}
    with the extra conservation law
    \begin{equation}
        \label{eq:Extra_conservation_law}
        \partial_t (\rho E) + \div F(\qq) = \div\left(\varepsilon\nabla\rho E\right).
    \end{equation}
    Here, the flux tensor is given by $\f(\qq) = (\rho \vv, \rho \vv \otimes \vv + p \mathbb{I}, S \vv)^T$, $F = (\rho E + p) \vv$ is the energy flux and $\PP = (0,\mathbf{0},\Pi)$ denotes a production source term acting only on the entropy.
    The Euler equations, as given above, are only one example in the general class of symmetric hyperbolic and thermodynamically compatible (SHTC) systems, see e.g. \cite{Godunov1961,Rom1998,Godunov2003}.
    Energy equation \eqref{eq:Extra_conservation_law} implies the compatibility condition
    \begin{equation}
        \label{eq:CompatibilitySHTC}
        \pp \cdot (\partial_t \qq+ \div \f(\qq)) = \partial_t (\rho E) + \div F(\qq)
    \end{equation}
    Note, that it holds $\pp \cdot \partial_t \qq= \partial_t{(\rho E)}$, where $\pp$ and $\qq$ are the thermodynamic dual variables and the state variables defined in \eqref{eq:PtoQ} and \eqref{eq:QtoP}.
    Thus, \eqref{eq:CompatibilitySHTC} reduces to the compatibility of the fluxes
    \begin{equation}
        \label{eq:CompatibilityFlux}
        \pp \cdot \left(\div \f(\qq)\right) = \div F(\qq).
    \end{equation}
    Consequently, we have the compatibility of the dissipative terms with the production term $\PP$
    \begin{equation}
        \label{eq:CompatibilityDissip}
        \pp \cdot \PP + \pp \cdot \div\left(\varepsilon\nabla \qq\right) = \div\left(\varepsilon\nabla\rho E\right).
    \end{equation}
    Following \cite{Abgrall2023HTC,Busto2022}, the compatible entropy production term $\Pi$ in the source term $\PP$ is given by
    \begin{equation}
        \label{eq:Entropy_production_continuous}
    { \Pi = \frac{\varepsilon}{T} \nabla \qq \cdot (\HH \nabla \qq) \geq 0.}
    \end{equation}
    For a convex total energy $\rho E$ and $\varepsilon \geq 0$ the production term is always non-negative since $\HH$ is symmetric and strictly positive definite (SPD).
    {We refer to \cite[Lemma 7.3 and Section 4.1.6]{Feireisl2021book} and \cite{LukArxiv2025} for the SPD property of $\HH$ .}
    Note that the Euler system \eqref{eq:Euler_in_S}, respectively \eqref{eq:Euler_compact}, contains a balance law for the entropy, since it is a state variable, whereas the conservation of the total energy \eqref{eq:Energy_conservation}, respectively \eqref{eq:Extra_conservation_law}, is a consequence.

    \section{The numerical scheme}\label{sec:NumericalScheme}
    \subsection{Spatial discretization}
    We use the notation forwarded in \cite{FeiLukMiz2019}.
    The space domain given by $\Omega_h = \Omega \subset \mathbb{R}^d, d = 1, 2, 3$ where $\Omega_h := [0, l]^d, l >0$ is divided into finite volume cells $K$, i.e.
    \begin{equation}
        \overline \Omega_h = \bigcup_{K\in \Th} \overline K.
    \end{equation}
    Therein $\Th$ is the set of all cells $K$ and forms a regular quadrilateral grid.
    In two space dimensions, the cell $K$ and its center $C_K$ corresponding to a uniform mesh size $h$ are given by
    \begin{align*}
        &K := \left[x_{i-1/2,j},x_{i+1/2,j}\right) \times \left[x_{i,j-1/2},x_{i,j+1/2}\right), \\ &
        C_K := (x_i, x_j) = \left(\frac{x_{i-1/2,j}+x_{i+1/2,j}}{{2}},\frac{x_{i,j-1/2}+x_{i,j+1/2}}{{2}}\right).
    \end{align*}

    Let the cell average on cell $K$ be expressed by the following projection
    \begin{equation}
        \label{eq:projection}
        P_h : L^1(\Omega) \to X(\Th), \quad (P_h(\varphi))_K := \frac{1}{h^d}\int_K \varphi(x)\ddd{x},
    \end{equation}
    where $X(\Th)$ denotes the space of piecewise constant functions defined on the mesh $\Th$.
    A piecewise function $g_h \in X(\Th)$ coincides on cell $K$ with the cell average $g_K = (P_h(g_h))_K$.
    Moreover we have
    \begin{equation}
        \label{eq:integral_on_Omega}
        \int_\Omega g_h \ddd{x} = h^d \sum_{K\in \Th} g_K.
    \end{equation}
    The common edge between two neighbouring cells $K$ and $L$ is denoted by {$\sigma=KL$} and the set of neighbouring cells is denoted by $\mathcal{N}(K)$.
    The set of all edges $\sigma$ of a cell $K \in \Th$ is denoted by $\mathcal{E}$.

    We define the following discrete operators with respect to cell $K$ applied to $g_h \in X(\Th)$ respectively.
    Since we consider a quadrilateral mesh, all operators are considered in direction $s$, where $s = 1, ..., d$.
    Centered differences in direction $s$ based on the neighbouring cells $L = K - h \ee_s$ and $R = K + h \ee_s$ where $\ee_s$ denotes the unit vector in the $s$-th space direction are then defined as
    \begin{equation}
        \left(\dhcentral g_h\right)_K = \frac{g_R - g_L}{2 h}
    \end{equation}
    and accordingly the respective upwind operators are given by
    \begin{equation}
        \left(\dhp g_h \right)_K= \frac{g_R - g_K}{h}, \quad \left(\dhm g_h \right)_K = \frac{g_K - g_L}{h}.
    \end{equation}
    Using the above given discrete operators, we define the discrete divergence operators of a piecewise vector function $\bm g_h = (g_h^1, \dots, g_h^d)^T$ based on cell averages by
    \begin{equation}
        \left(\divh \bm g_h\right)_K = \sum_{s=1}^{d} \left(\dhcentral g_h^s\right). 
    \end{equation}
    Further we define a jump over an interface $\sigma = LK$, with $\nn_{LK}$ being the normal outward pointing vector from {cell $L$ to $K$} by
    \begin{equation}
        \jump{g_h}_\sigma := g_K \nn_{LK} + g_L \nn_{KL} = (g_K - g_L)\nn_{LK}.
    \end{equation}
    The mean value over the interface $\sigma = LK$ reads
    \begin{equation}
        \mean{g_h}_\sigma := \frac{g_L + g_K}{2},
    \end{equation}
    and the mean value on cell $K$
    \begin{equation}
        \mean{g_h}_K^s := \frac{g_R + g_L}{2}, \quad L = K - h \ee_s, \quad R = K + h \ee_s.
    \end{equation}
    \subsection{The hyperbolic thermodynamically compatible FV scheme}
    Our aim is to prove the convergence of the novel hyperbolic thermodynamically compatible finite volume (HTC-FV) scheme given in \cite{Abgrall2023HTC,Thomann2023}.
    Throughout this work, we will consider the space semi-discrete version of the HTC-FV scheme, thus the time remains continuous.
    The space-semi discrete scheme can be equipped with an explicit high order time discretization based on the method of lines approach or a fully discrete implicit scheme can be derived, see \cite{Abgrall2023HTC,Busto2022,Thomann2023} for details.

    The numerical flux for the HTC-FV scheme at the interface $\sigma = LK$, where $L = K - h \ee_s$, is given by
    \begin{equation}
        \label{eq:HTC_num_flux}
        \f_{h,\sigma} = \mean{\f^s(\qq_h)}_\sigma - \alpha_\sigma \jump{\pp_h}_\sigma - \lambda_\sigma \jump{\qq_h}_\sigma,
    \end{equation}
    where
    \begin{equation}
        \label{eq:lambda_sig}
    \lambda_\sigma := \frac{1}{2}\max_{k}\left(|\lambda_k(\qq_L)|,|\lambda_k(\qq_K)|\right)
   \end{equation}
    denotes the local diffusion coefficient independent of $h$.
    We want to emphasize, that this numerical flux consists of the \textit{novel Abgrall flux} and a Rusanov-type dissipation.
    In \eqref{eq:HTC_num_flux} $\lambda_k$ denotes the $k-th$ eigenvalue of the flux Jacobian $\nabla_{\qq} \f^s$ in direction $s$ given in \eqref{eq:Eigenvalues}.
    Instead of a local diffusion coefficient, a global diffusion coefficient can be chosen, e.g. as in the Lax-Friedrichs scheme
    \begin{equation}\label{eq:lambda}
        \lambda := \frac{1}{2}\max_{K \in \Th}\max_{k}\left(|\lambda_k(\qq_L)|,|\lambda_k(\qq_K)|\right).
    \end{equation}
    The diffusion terms $\lambda_\sigma\jump{\qq_h}_\sigma$ correspond to a consistent discretization of the vanishing viscosity regularization terms introduced in \eqref{eq:Euler_in_S_vvr}.
    It stabilizes the centred term $\mean{\f^s(\qq_h)}_\sigma$ which will yield a consistent discretization of the divergence of the flux function $\f$ in \eqref{eq:Euler_in_S_vvr}.
    Finally, the so called \textit{Abgrall correction term} $\alpha_\sigma\jump{\pp_h}_\sigma$ ensures the discrete flux compatibility corresponding to \eqref{eq:CompatibilityFlux}, see \cite{Abgrall2023HTC} for details.
    Therein, the so called \textit{Abgrall correction factor} $\alpha_\sigma$ evaluated at the interface $\sigma = LK$ is given by
    \begin{equation}
        \label{eq:AbgrallCorrectionFactor}
        \alpha_\sigma^s = \frac{\jump{F^s(\qq_h)}_\sigma + \mean{\f^s(\qq_h)}_\sigma\cdot \jump{\pp_h}_\sigma - \jump{\pp_h \cdot \f^s(\qq_h)}_\sigma}{\jump{\pp_h}_\sigma \cdot \jump{\pp_h}_\sigma},
    \end{equation}
    where $F^s$ denotes the energy flux in direction $s$ of the interface, i.e. $F^s(\qq_h) := \frac{\m_h\cdot \nn_\sigma}{\rho_h}\left((\rho E)_h + p_h\right)$.
    Analogous notation holds for $\f^s$.

    Defining a Roe-type matrix $(\HH_h)_\sigma$ allows to rewrite jumps in state variables $\jump{\qq_h}_\sigma$ over an interface $\sigma = LK$  in terms of dual variables $\jump{\pp_h}_\sigma$ as follows
    \begin{equation}
        \label{eq:Roe_property}
        \left({\HH}_h\right)_\sigma \jump{\qq_h}_\sigma = \jump{\pp_h}_\sigma .
    \end{equation}
    Following \cite{Abgrall2023HTC,Busto2022}, the Roe matrix $(\HH_h)_\sigma$ is constructed by defining a simple segment path $\mathbf{\psi}$ in $\qq$ given by
    \begin{equation}
        \mathbf{\psi}(s) = \qq_L + s(\qq_K - \qq_L), \quad 0\leq s \leq 1
    \end{equation}
    which yields by integration the Roe matrix
    \begin{equation}
        \label{eq:Roe_matrix}
        (\HH_h)_\sigma = \int_0^1 \HH(\mathbf{\psi}(s)) \ddd{s}.
    \end{equation}
    Using \eqref{eq:Roe_property}, we can rewrite the numerical flux \eqref{eq:HTC_num_flux} as follows
    \begin{equation}
        \label{eq:HTC_num_flux_q}
        \f_{h,\sigma} = \mean{\f^s(\qq_h)}_\sigma - \left(\alpha_\sigma \left({\HH}_h\right)_\sigma + \lambda_\sigma \mathbb{I}\right) \jump{\qq_h}_\sigma.
    \end{equation}
    Note that $\alpha_\sigma$ is not contributing to the numerical dissipation in general, since it does not have a fixed sign.
    Analogously, the correction factor \eqref{eq:AbgrallCorrectionFactor} can be rewritten in jumps in $\qq$ as follows
        \begin{equation}
        \label{eq:AbgrallCorrectionFactor_q}
        \alpha_\sigma^s = \frac{\jump{F^s(\qq_h)}_\sigma + \mean{\f^s(\qq_h)}_\sigma\cdot \left({\HH}_h\right)_\sigma \jump{\qq_h}_\sigma - \jump{\pp_h \cdot \f^s(\qq_h)}_\sigma}{\left({\HH}_h\right)_\sigma \jump{\qq_h}_\sigma \cdot \left({\HH}_h\right)_\sigma \jump{\qq_h}_\sigma}.
    \end{equation}

    Finally, the space-semi discrete HTC-FV scheme for the Euler system \eqref{eq:Euler_in_S} is given for $K \in \Th$ and $t>0$ by
    \begin{subequations}
        \label{eq:HTC_scheme}
    \begin{eqnarray}
        &&\frac{\partial \rho_K(t) }{\partial t} + \left(\divh \mm_h\right)_K - \frac{1}{h} \sum_{\sigma \in \pK}\left(\alpha_\sigma \jump{r_h(t)}_\sigma+\lambda_\sigma \jump{\rho_h(t)}_\sigma\right) = 0 \label{disc_continuity}\\
        &&
        \begin{split}
        \frac{\partial \mm_K(t) }{\partial t} &+ \left(\divh \left(\frac{\mm_h(t) \otimes \mm_h(t)}{\rho_h(t)}+ p_h(t)\mathbb{I}\right)\right)_K \\
        &- \frac{1}{h} \sum_{\sigma \in \pK}\left(\alpha_\sigma\jump{\vv_h(t)}_\sigma + \lambda_\sigma \jump{\mm_h(t)}_\sigma\right) = 0
        \end{split}\label{disc_momentum}\\
        &&
        \begin{split}\frac{\partial S_K(t) }{\partial t} &+ \left(\divh \left(\frac{\mm_h(t) S_h(t)}{\rho_h(t)}\right)\right)_K \\
            &- \frac{1}{h} \sum_{\sigma \in \pK}\left(\alpha_\sigma \jump{T_h(t)}_\sigma + \lambda_\sigma \jump{S_h(t)}_\sigma\right) = \frac{1}{h}(\Pi_h)_K.
        \end{split} \label{disc_entropy}
    \end{eqnarray}
    \end{subequations}
    A consistent discretization of the entropy production term $\Pi$ which is compatible with the compatibility condition on the dissipative processes \eqref{eq:CompatibilityDissip} is given by
    \begin{equation}
        \label{eq:Entropy_prod_K}
        (\Pi_h)_K = \frac{1}{2} \frac{1}{T_K}\sum_{\sigma \in \pK} \lambda_\sigma \jump{\qq_h(t)}_\sigma {\cdot}\left(\HH_h(t)\right)_\sigma \jump{\qq_h(t)}_\sigma.
    \end{equation}
    Note, that $\varepsilon$ in the vanishing viscosity formulation is discretized at the interface and due to \eqref{eq:HTC_num_flux} this corresponds to setting $\varepsilon_\sigma = h\lambda_\sigma$.
    By construction solutions $\qq_h = (\rho_h, \mm_h, S_h)$ of the numerical scheme \eqref{eq:HTC_scheme} fulfil the additional compatible discrete energy conservation law
    \begin{equation}
        \label{eq:HTC_energy}
        \frac{\partial (\rho E)_K(t)}{\partial t} + (\divh F(\qq_h))_K ~-~ \frac{1}{h} \sum_{\sigma \in \pK} \lambda_\sigma \jump{\rho E(t)}_\sigma = 0.
    \end{equation}
    For details on the discrete compatibility of \eqref{eq:HTC_scheme} using the entropy production term discretization \eqref{eq:Entropy_prod_K} with the discrete energy conservation law \eqref{eq:HTC_energy} see the proofs given in \cite{Busto2022,Abgrall2023HTC,Thomann2023}.
    \begin{remark}
        Even though the jump term associated with $\alpha_\sigma$ is written together with the numerical diffusion in \eqref{eq:HTC_num_flux_q}, we want to emphasize that $\alpha_\sigma$ does not have a sign and therefore does not contribute to the diffusivity of the numerical scheme in general.
        However, changing the definition of the numerical diffusion coefficient $\lambda_\sigma$ the positivity of the term $\alpha_\sigma \left({\HH}_h\right)_\sigma + \lambda_\sigma$ in the numerical flux can be ensured. This is subject of the following lemma.
        Note that the compatibility conditions require the analogue modification of the dissipation in the entropy production term, when $\lambda_\sigma$ in the numerical flux is changed.
    \end{remark}
        \begin{lemma}
        \label{lem:Positivity_of_diffusive_terms}
        Let $T_K (t) > 0$ uniformly and $(\HH_h)_\sigma$ given by \eqref{eq:Roe_matrix}.
        Consider the modified numerical flux \eqref{eq:HTC_num_flux_q} given by
        \begin{equation}
           \f_{h,\sigma} = \mean{\f^s(\qq_h)}_\sigma - \left(\alpha_\sigma \left({\HH}_h\right)_\sigma + \left(\lambda_\sigma + |\min(0,\alpha_\sigma)| \overline{\kappa}\right)\mathbb{I}\right) \jump{\qq_h}_\sigma,
        \end{equation}
        where the constant $\overline{\kappa}$ being an upper bound on the eigenvalues of $~\HH$.

        Then the coefficient $\alpha_\sigma \left({\HH}_h\right)_\sigma + \left( |\alpha_\sigma| \overline{\kappa} + \lambda_\sigma\right) \mathbb{I}$ in \eqref{eq:HTC_num_flux_q} is positive.
    \end{lemma}
    \begin{proof}
        The result follows immediately if $\alpha_\sigma \geq 0$. Thus, let $\alpha_\sigma < 0$.
        Since $\HH$ is {SPD}, there exists a matrix $\mathbf{V}$ and a diagonal matrix $\mathbf{K}$ such that $\HH = \mathbf{V} \mathbf{K} \mathbf{V}^{-1}$.
        Let $\overline{\kappa}$ be constant with $\overline{\kappa} \geq \max_{k}\max(\kappa_k(\qq_L),\kappa_k(\qq_R))$ where $\kappa_k(\qq)$ denote the entries of the diagonal matrix $\mathbf{K}$.
        Then $(\HH_h)_\sigma \leq \overline{\kappa}\mathbb{I}$ and it follows
        \begin{equation}
            \alpha_\sigma (\HH_h)_\sigma + \left(\lambda_\sigma +|\alpha_\sigma|\kappa\right)\mathbb{I} \geq \lambda_\sigma \geq 0
        \end{equation}
        which concludes the proof.
    \end{proof}
    We conclude this section with a result on the Abgrall correction factor $\alpha_\sigma$.
    \begin{lemma}
    \label{lem:Abgrall_alpha}
    The {Abgrall} correction factor {$\alpha_\sigma^s$ defined in \eqref{eq:AbgrallCorrectionFactor} over an interface $\sigma$ in direction $s = 1, \dots, d$ converges as $\jump{\qq_h}_\sigma\to 0$.}
    \end{lemma}

    \begin{proof}
    Let us consider without loss of generality the interface $\sigma = LK$.
    {
    We expand each term in the nominator of \eqref{eq:AbgrallCorrectionFactor} by means of the Taylor series as follows.
    The first term regarding the difference of the total energy flux can be estimated by
        \begin{equation}
        \jump{F^s(\qq_h)}_\sigma = \nabla_{\qq} F^s(\qq_L) \cdot \jump{\qq_h}_\sigma + \frac{1}{2} \nabla_{\qq}^2 F^s(\qq_L) \jump{\qq_h}_\sigma\cdot  \jump{\qq_h}_\sigma + \mathcal{O}(\jump{\qq_h}_\sigma^3).
        \end{equation}
    For the second term we obtain
        \begin{equation}
        \mean{\f^s(\qq_h)}_\sigma\cdot \jump{\pp_h}_\sigma = \left(\f^s(\qq_L) + \nabla_{\qq} \f^s(\qq_L) \jump{\qq_h}_\sigma  \right) \cdot \jump{\pp_h}_\sigma+ \mathcal{O}(\jump{\qq_h}_\sigma^3).
        \end{equation}
    The last term in the nominator yields the following expansion
    \begin{multline}
         -\jump{\pp_h \cdot \f^s(\qq_h)}_\sigma = -\jump{\pp_h}_\sigma \cdot \f^s(\qq_R) - \pp_L \cdot \left( \nabla_{\qq} \f^s(\qq_L) \jump{\qq_h}_\sigma + \frac{1}{2}D^2 \f^s(\qq_L) \cdot \jump{\qq_h}_\sigma \otimes \jump{\qq_h}_\sigma\right) + \mathcal{O}(\jump{\qq_h}_\sigma^3),
    \end{multline}
    where $D^2 \f^s(\qq_L) = \frac{\partial^2\f^s_i}{\partial \qq_j\partial\qq_k} e_i \otimes e_j \otimes e_k$ denotes the third order tensor containing the second order derivatives of $\f^s$ with respect to $\qq$.
    Gathering the first order terms in $\jump{\pp_h}_\sigma$ we obtain using again Taylor expansion
    \begin{equation*}
        -\jump{\pp_h}_\sigma \cdot (\f^s(\qq_R) - \f^s(\qq_L)) = -\jump{\pp_h}_\sigma \cdot \nabla_{\qq} \f^s(\qq_L) \jump{\qq_h}_\sigma + \mathcal{O}(\jump{\qq_h}_\sigma^3).
    \end{equation*}
    For the first order terms in $\jump{\qq_h}_\sigma$ vanish due to the compatibility condition \eqref{eq:CompatibilityFlux}.
    Making use of the transform from $\pp$ to $\qq$ variables and that  and applying the rule of L'Hospital, we obtain
        \begin{equation}
         \lim\limits_{\jump{\qq_h}_\sigma\to 0} \alpha_\sigma^s = \lim\limits_{\jump{\qq_h}_\sigma\to 0}  \frac{\frac{1}{2} \nabla_{\qq}^2 F^s(\qq_L) \ee_s \cdot \ee_s +   \nabla_{\qq} \f^s(\qq_L)\ee_s \cdot (\HH_h)_\sigma\ee_s - \pp_L \cdot  \frac{1}{2}D^2 \f^s(\qq_L) \cdot \ee_s \otimes \ee_s}{(\HH_h)_\sigma \ee_s \cdot (\HH_h)_\sigma \ee_s}.
        \end{equation}
    Note that it is equivalent to consider the limit $\jump{\pp_h}_\sigma \to 0$.
    Denoting by $(\hat\HH_h)_\sigma$ the limit of $(\HH_h)_\sigma$ as $\jump{\qq_h}_\sigma \to 0$, the limit of $\alpha_\sigma^s$ is given by
    \begin{equation*}
    \frac{\frac{1}{2} \nabla_{\qq}^2 F^s(\qq_L) \ee_s \cdot \ee_s +   \nabla_{\qq} \f^s(\qq_L)\ee_s \cdot (\hat \HH_h)_\sigma\ee_s - \pp_L \cdot  \frac{1}{2}D^2 \f^s(\qq_L) \cdot \ee_s \otimes \ee_s}{(\hat\HH_h)_\sigma \ee_s \cdot (\hat\HH_h)_\sigma \ee_s}.
    \end{equation*}
}
    \end{proof}

    \subsection{Upper and lower bounds of conserved and thermodynamical variables}

    In this section, we will derive upper and lower bounds for the state variables $\qq$ and the associated thermodynamical variables under the following assumption
    \begin{assumption}\label{ass:uniform_positivity}
We assume that there exist $\underline{\rho} > 0,\, \overline{E} > 0$ such that
        \begin{equation}
            0 < \underline{\rho} \leq \rho_h, \ \mbox{ and } \ (\rho E)_h \leq \overline{E}
        \end{equation}
        uniformly for $h \to 0, t \in [0,T]$.
    \end{assumption}
    Under this assumption we prove the non-negativity of the entropy production term.
     \begin{lemma}[Estimate of the entropy production term]
        \label{lem:Positivity_of_entropy_production}
        Let $(\HH_h)_\sigma$ given by \eqref{eq:Roe_matrix}.

        Then under Assumption \ref{ass:uniform_positivity}, the entropy source term defined in \eqref{eq:Entropy_prod_K} is non-negative for all $K \in \Th$.
    \end{lemma}
    \begin{proof}
    Since the energy Hessian is {SPD}, we can estimate
        \begin{equation}
        \jump{\qq_h}_\sigma {\cdot} \jump{\pp_h}_\sigma =  \jump{\qq_h}_\sigma \cdot {(\HH_h)_\sigma} \jump{\qq_h}_\sigma \geq 0.
        \end{equation}
    Due to the definition of temperature  in \eqref{eq:pressure_EOS}, and the uniform lower bound on the density, {$T_K$ for all $K \in \Th$ is positive
        \begin{equation}
        e(\rho_h,s_h) = \frac{\rho_h^{\gamma-1}}{\gamma - 1} \exp\left(\frac{s_h}{c_v}\right) > 0.
        \end{equation}
        Together with the positivity of $\lambda_\sigma$ we obtain $(\Pi_h)_K \geq 0$ for all $K \in \Th$.}
    \end{proof}

    \begin{lemma}
    \label{lem:3}
     {Under Assumption \ref{ass:uniform_positivity} it holds
        \begin{equation}
             0 < \underline{\rho} \leq \rho_h(t) \leq \overline{\rho}, \quad  |\m_h(t)| \leq \overline{m}, \quad 0 < \underline{p} \leq p_h(t)  \leq \overline{p}, \quad \underline{S} \leq S_h(t) \leq \overline{S}
        \end{equation}
         uniformly for $h \to 0$, $t \in [0,T]$ with constants $\overline{\rho}, \overline{m}, \underline{T}, \overline{T}, \underline{p}, \overline{p}, \underline{S}, \overline{S}$.
        }
    \end{lemma}
    \begin{proof}
     From the equation of state \eqref{eq:pressure_EOS}, and using the upper bound for the total energy, we have
        \begin{equation}
             p_h = (\gamma - 1) \left((\rho E)_h - \frac{1}{2}\frac{{|}\mm_h{|}^2}{\rho_h}\right) \leq (\gamma - 1) (\rho E)_h \leq (\gamma - 1) \overline{E} =: \overline{p}
        \end{equation}
    and
    \begin{equation}\label{eq:entropy_trans}
         T_h = \frac{1}{c_v \rho_h} \left((\rho E)_h - \frac{1}{2}\frac{{|}\mm_h{|}^2}{\rho_h}\right) \leq \frac{1}{c_v \underline{\rho}}\overline{E} =: \overline{T}.
        \end{equation}
    {To show the minimal entropy principle, we consider the renormalized entropy inequality.
        It is derived from the transport equation for the specific entropy $s$
        \begin{equation}
            \partial_t s  + \vv \cdot \nabla s = \nabla \cdot (\varepsilon \nabla s) + 2 \frac{\varepsilon}{\rho} \sum_{k=1}^d \frac{\partial s}{\partial x_k}\frac{\partial \rho}{\partial x_k} + \frac{1}{\rho} \Pi
        \end{equation}
        which is obtained by combining \eqref{eq:Euler_in_S_vvr:dens} and \eqref{eq:Euler_in_S_vvr:entropy}.
         Consider now the specific normalized entropy $s_\chi = \chi \circ s$ introduced in \cite{Harten1983} where $\chi$ is an increasing bounded concave function, i.e.
        \begin{equation}
        \chi^\prime(s) \geq 0, \quad \chi^{\prime\prime}(s) \leq 0.
        \end{equation}
        Multiplying \eqref{eq:entropy_trans} by $\rho \chi^\prime(s)$ and adding the density equation \eqref{eq:Euler_in_S_vvr:dens} multiplied by $s_\chi$ we obtain the renormalized entropy inequality for $S_\chi := \rho s_\chi$ given by
        \begin{equation}\label{eq:entropy_renorm}
            \partial_t S_\chi  + \nabla \cdot (S_\chi \vv ) = \nabla \cdot (\varepsilon \nabla S_\chi) + \chi^\prime \Pi - \varepsilon \rho \chi^{\prime\prime} \sum_{k=1}^d \left(\frac{\partial s}{\partial x_k}\right)^2.
        \end{equation}
        Since $\chi^\prime \geq 0$ and $\chi^{\prime\prime} \leq 0$ the second and third term on the right-hand-side are non-negative.
        Moreover, the first term on the right-hand-side corresponds to the numerical diffusion in the numerical flux and thus the renormalized entropy inequality is achieved.
        \\
        In \cite[Theorem 1.2]{Lax1971} it was shown by Lax that the (local) Lax-Friedrichs scheme for all entropies satisfies the discrete entropy inequality, thus also for $S_\chi$.
        Since with the choice of dissipation from Lemma \ref{lem:Positivity_of_diffusive_terms} the HTC scheme can be written as the local Lax-Friedrichs scheme. Similarily, as in \cite[Theorem 1.2]{Lax1971}, it then satisfies also the discrete renormalized entropy inequality.
        \\
        To obtain the minimum entropy principle, we follow the argumentation of Lax \cite[Section 3]{Lax1971} and Tadmor \cite{Tadmor1986} by considering
        \begin{equation}
            \chi(s) = \min(s - s_0, 0)
        \end{equation}
        with $s_0:= \min_K s_h(0), \overline{C} = \max_K \frac{\rho_h(0) T_0}{T^{1/(\gamma - 1)}}$ and $T_0 = \frac{1}{(\gamma - 1) c_v}$.
        For more details on this approach see also \cite{FeiLukMiz2019,Brezina2018,Lukacova2023}.
       \\
        Due to $\chi^\prime \geq 0$ and $\chi^{\prime\prime} \leq 0$ the positivity of the source terms is ensured and by integration over $\Omega$ and $(0,T)$ we obtain
        \begin{equation}
           \int_\Omega \rho_h \chi(s_h) d\x \geq \int_\Omega \rho(0)_h \chi(s(0)_h) d\x \geq \underline{\rho}\left| \Omega\right| \chi(s_0),
        \end{equation}
        which yields the minimum entropy principle
        \begin{equation}
            s_h(t) \geq s_0 = \min s_h(0).
        \end{equation}
        Further we have
    \begin{align}
            \begin{split}
            &c_v (\gamma - 1)  \log\left(\frac{\left(\frac{T_h(t)}{T_0}\right)^{1/(\gamma - 1)}}{\rho_h(t)}\right) =  c_v \log\left(\frac{ T_h(t)}{T_0\rho_h(t)^{\gamma-1}}\right) =  s_h(t)  \\
            & \geq s_0 = \min_K \left(c_v (\gamma - 1)  \log\left(\frac{\left(\frac{T_h(0)}{T_0}\right)^{1/(\gamma - 1)}}{\rho_h(0)}\right) \right) \geq c_v (\gamma - 1) \log \left(\frac{1}{\overline{C}}\right).
            \end{split}
    \end{align}
        Thus, we obtain the estimate for the density
        \begin{equation}
            0 < \underline{\rho} \leq \rho_h(t) \leq \overline{C} T_0^{-1} \overline{T}^{1/(\gamma -1)} = \overline{\rho}.
        \end{equation}
        Using this result, we obtain the upper bound for the momentum}
        \begin{equation}
               {|\mm_h|}^2 = 2 \rho_h\left((\rho E)_h - (\rho e)_h\right) \leq 2 \rho_h (\rho E)_h \leq 2{\overline{\rho}}\overline{E} = \overline{m}.
        \end{equation}
        {Furthermore, we have via the equation of state
        \begin{equation}
            T_h(t) = \frac{1}{c_v(\gamma - 1)}\rho_h(t)^{\gamma-1} \exp(s_h(t)/c_v) \geq \frac{1}{c_v(\gamma - 1)}\underline{\rho}^{\gamma-1} \exp(s_0/c_v) = \underline{T}
        \end{equation}
        and
        \begin{equation}
            p_h(t) = (\gamma - 1) c_v \rho_h(t) T_h(t) \geq (\gamma - 1) c_v \underline{\rho}\underline{T} = \underline{p}.
        \end{equation}
        Finally we can set $\underline{S}: = \underline{\rho}\underline{s} \leq S_h(t) \leq \overline{\rho}\overline{s} = :\overline{S}$ which concludes the proof.
    }
    \end{proof}

    \subsection{Stability of the Numerical Scheme}\label{sec:stability}
    In this section, we summarize the stability of the numerical scheme \eqref{eq:HTC_scheme} with \eqref{eq:HTC_energy} by stating a priori estimates.
    Summing up the continuity equation \eqref{disc_continuity} over all $K \in \Th$ and multiplying by $h^d$, we obtain with \eqref{eq:integral_on_Omega} and periodic boundary conditions
    \begin{equation}
        \label{eq:mass_conservation_on_omega}
        \int_\Omega \rho_h(t) \ddd{x} = \int_\Omega \rho_h(0) \ddd{x}.
    \end{equation}
    Analogously, integrating \eqref{eq:HTC_energy} over $\Omega$ we obtain the discrete energy conservation
      \begin{equation}
        \label{eq:energy_conservation_on_omega}
        \int_\Omega (\rho E)_h(t) \ddd{x} = \int_\Omega (\rho E)_h(0) \ddd{x}.
    \end{equation}
    In the previous section, we have shown that the density and temperature (pressure) are uniformly bounded from below away from zero, thus the same applies for the total energy, i.e. $(\rho E)_h(t) > 0$ for all $t \in [0,T]$.
    With the positivity of $(\Pi_h)_K$, see Lemma \ref{lem:Positivity_of_entropy_production} combined with Lemma \ref{lem:3}, we obtain for the entropy
    \begin{equation}
        \label{eq:entropy_production_on_omega}
        \int_\Omega S_h(t) \ddd{x} \geq \int_\Omega S_h(0) \ddd{x}.
    \end{equation}

    The conservation of mass \eqref{eq:mass_conservation_on_omega} and total energy \eqref{eq:energy_conservation_on_omega} and an upper bound on the density and total energy from Lemma \ref{lem:3} give the following a priori estimates
    \begin{align}
        \label{eq:apriori_estimates}
        \begin{split}
        \rho_h &\in L^\infty((0,T)\times\Omega), \quad
        p_h \in L^\infty((0,T)\times\Omega), \quad (\rho E)_h \in L^\infty((0,T)\times\Omega)\\
        S_h &\in L^\infty((0,T)\times\Omega), \quad \sqrt{\rho_h}\vv_h \in L^\infty((0,T)\times\Omega), \quad \m_h = \rho_h \vv_h \in L^\infty((0,T)\times\Omega).
        \end{split}
    \end{align}

    \section{Consistency}\label{sec:Consistency}
    The aim of this section is to prove the consistency of the HTC scheme \eqref{eq:HTC_scheme}.
    To establish the consistency errors, we multiply the continuity and entropy equation given by \eqref{disc_continuity} and \eqref{disc_entropy} by $h^d (P_h \varphi(t))_K$ where $\varphi \in C^1([0,T],C^2(\Omega))$ and the momentum equation \eqref{disc_momentum} by $h^d (P_h \phibf(t))_K$ with $\phibf \in C^1([0,T], C^2(\overline \Omega; \mathbb{R}^d))$.
    For the sake of a compact notation, we consider with a slight abuse of notation $\phibf \in C^1([0,T], C^2(\Omega; \mathbb{R}^{d+2}))$ for all variables $\qq$.
    Note that $P_h$ denotes the projection operator defined in \eqref{eq:projection}.
    Then we sum the resulting equations over $K \in \Th$ and integrate in time.
    We will analyse, following \cite{FeiLukMiz2019}, the terms associated to time derivatives, convection and numerical diffusion.
    Due to the derivation of the HTC scheme, we have to consider in addition the flux correction terms and the entropy production term which does not appear in the previously studied entropy stable FV schemes \cite{FeiLukMiz2019}.

    System \eqref{eq:Euler_compact} is usually equipped with suitable boundary conditions. Here, we will only consider periodic boundary conditions, as done in \cite{FeiLukMiz2019}.
    Thus the computational domain is given by the flat torus
    \begin{equation}
        \label{eq:Periodic_domain}
        \Omega = \left([0,1]|_{\left\lbrace0,1\right\rbrace}\right)^d.
    \end{equation}
    Further, we prescribe regular initial data given by
    \begin{equation}
        \label{eq:init_data}
        \rho^0 \in L^{\infty}(\Omega), \quad \rho^0 > 0, \quad \m^0 \in L^\infty(\Omega; \mathbb{R}^d), \quad S^0 \in L^\infty(\Omega),
    \end{equation}
    such that
    \begin{equation}
        \label{eq:init_data_thermodyn}
        p^0 = \left(\rho^0\right)^\gamma \exp\left(\frac{S^0}{c_v \rho^0}\right) > 0 \quad\text{and}\quad T^0 = \frac{(\rho^0)^{\gamma - 1}}{c_v(\gamma - 1)} \exp\left(\frac{S^0}{c_v \rho^0}\right) > 0.
    \end{equation}

    \subsection{Consistency errors}

    Regarding the time derivatives using integration by parts we obtain
    \begin{equation}
        h^d \int_0^T \frac{\ddd{}}{\ddd{t}}\sum_{K\in \Th} \qq_K(t) \left(P_h\phibf(t)\right)_K \ddd{t} = \left[\int_\Omega \qq_h(\tau)\cdot \phibf(\tau,x)\ddd{x}\right]_{t = 0}^{t=T} - \int_0^T \int_\Omega \qq_h \cdot \partial_t \phibf(t) \ddd{x}\ddd{t}.
    \end{equation}
    Next, we consider the convective terms.
Performing the Taylor series expansion {component by component} for $\phibf$, we obtain
    \begin{equation}
        h^d \int_0^T \sum_{K\in \Th} \left(\divh \f(\qq_h(t))\right)_K \left(P_h\phibf(t)\right)_K \ddd{t} = - \int_0^T \int_\Omega \f(\qq_h ): \nabla \phibf(t,x) \ddd{x}\ddd{t} + \rr,
    \end{equation}
    where $\rr = (r_1, ..., r_{d+1})$ are bounded componentwise by
    \begin{align}
        &r_1 \leq h \left\|\frac{\ddd{^2\varphi}(\hat{x})}{\ddd{x}^2}\right\|_{C(0,T)} \|\m\|_{L^\infty}, \\
        &r_k \leq h \left\|\frac{\ddd{^2\phibf}(\hat{x})}{\ddd{x}^2}\right\|_{C(0,T)} \left(\|\sqrt{\rho_h(t)}\vv_h(t)\|_{L^\infty} + \|p_h(t)\|_{L^\infty}\right), \quad k = 2, \dots d+1 \\
        &r_{d+2} \leq h \left\|\frac{\ddd{^2\varphi}(\hat{x})}{\ddd{x}^2}\right\|_{C(0,T)} \|S_h(t) \vv_h(t)\|_{L^\infty}.
    \end{align}
    Therein, we have used the following notation
    \begin{equation}
        \frac{\ddd{^2\varphi}(\hat{x})}{\ddd{x}^2} = \left(\frac{\partial^2 \varphi}{\partial x_i \partial x_j}\right)_{i,j = 1}^d(\hat{x})
    \end{equation}
    with $\hat{x}\in (x-h\ee_s, x + h\ee_s)$  from the Taylor series expansion.

    The consistency of the numerical diffusion terms
    \begin{equation}
        \label{eq:diffusion_consistency}
        h^{d-1} \int_0^T \sum_{K\in \Th}\sum_{\sigma \in \pK} \lambda_\sigma \jump{\qq_h(t)}_\sigma (P_h \phibf(t))_K \ddd{t}
    \end{equation}
    follows  analogously as in \cite[Section~9.3]{FeiLukMiz2019} provided that the weak bounded variation (BV) property holds which we will prove in the following.
    \begin{definition}[weak BV condition]
        Let $\qq_h(t)$ be a solution of the HTC scheme \eqref{eq:HTC_scheme}.
        Then $\qq_h(t)$ is weakly BV if
        \begin{equation}
            \label{eq:weakBV}
            \int_0^T \sum_{\sigma \in \mathcal{E}} \lambda_\sigma \left|\jump{\qq_h(t)}_\sigma\right| h^d \to 0
        \end{equation}
        as $h \to 0^+$, where $\lambda_\sigma$ denotes the numerical viscosity coefficient given in \eqref{eq:lambda_sig}.
    \end{definition}
    \begin{lemma}[weak BV property]
        \label{lem:weakBV}
        Let $\qq_h(t)$ be a solution of the HTC scheme \eqref{eq:HTC_scheme} and let Assumption \ref{ass:uniform_positivity} hold.

        Then the weak BV condition \eqref{eq:weakBV} holds.
    \end{lemma}
    \begin{proof}
    {By Lemma  \ref{lem:Positivity_of_entropy_production}, we know that the entropy source term $\Pi_h$ is non-negative.}
        Thus integrating the entropy equation \eqref{disc_entropy} over $\Omega$ with periodic boundary conditions and $(0,T)$ yields
        \begin{equation}
            \int_\Omega S_h(T)\ddd{x} - \int_0^T \sum_{\sigma \in \mathcal{E}}h^{d-1}(\Pi_h(t))_\sigma \ddd{t} = \int_\Omega S_h(0)\ddd{x}.
        \end{equation}
    Due to Lemma \ref{lem:3} $S_h(\tau)$ for all $\tau \in [0,T]$ is bounded below by a constant $\underline{S}$ and from above by some constant $\overline{S}$.
        Thus, passing to the limit it follows
        \begin{equation}
            \label{eq:h_limit_Pi}
            - \int_0^T \sum_{\sigma \in \mathcal{E}}h^{d}(\Pi_h(t))_\sigma \ddd{t} \to 0^- \quad \text{for} \quad h \to 0^+.
        \end{equation}
    Further, due  $\HH$ {being SPD}, there exists {a uniform}  $\eta > 0$ such that $ \eta \jump{\qq_h}_\sigma\leq(\HH_h)_\sigma \jump{\qq_h}_\sigma = \jump{\pp_h}.$
    {This can be seen by taking the minimum over all eigenvalues of ${(\HH_h)_\sigma}$ which are uniformly bounded away from zero due to Lemma \ref{lem:3}.}
    Consequently, it holds
        \begin{equation}
            \label{eq:Estimate_Pi}
            \int_0^T \sum_{\sigma \in \mathcal{E}} h^d \frac{\lambda_\sigma}{2} \frac{1}{T_K} \jump{\qq_h(t)}_\sigma \jump{\pp_h(t)}_\sigma \ddd{t} \geq \frac{\eta}{2}\int_0^T \sum_{\sigma \in \mathcal{E}} h^d \lambda_\sigma \frac{1}{T_K} \jump{\qq_h(t)}_\sigma \jump{\qq_h(t)}_\sigma \ddd{t} \geq 0,
        \end{equation}
        where the first term tends to zero due to \eqref{eq:h_limit_Pi}.
        To show the BV estimate \eqref{eq:weakBV}, we apply the H\"older inequality
        \begin{equation}
            \label{eq:Aux3}
            \int_0^T \sum_{\sigma \in \mathcal{E}} \lambda_\sigma \Big|\jump{\qq_h(t)}_\sigma\Big| h^d \ddd{t}
            \leq \left(\int_0^T  \sum_{\sigma \in \mathcal{E}} h^d T_K \lambda_\sigma\ddd{t}\right)^{1/2} \left(\int_0^T\sum_{\sigma \in \mathcal{E}}h^d\lambda_\sigma \frac{1}{T_K}\Big|\jump{\qq_h(t)}\Big|^2\ddd{t}\right)^{1/2}.
        \end{equation}
        The second term on the right hand side of the inequality tends to zero due to \eqref{eq:h_limit_Pi} and \eqref{eq:Estimate_Pi}.
        Using the discrete trace inequality on the integral in the first term on the right hand side yields
        \begin{equation}
        \begin{split}
        \int_0^T  \sum_{\sigma \in \mathcal{E}} h^d T_K\lambda_\sigma\ddd{t} &\leq h \int_0^T \sum_{K\in \Th}\sum_{\sigma \in \pK} \int_\sigma T_K\lambda_\sigma \ddd{S} \ddd{t}\leq h \int_0^T \sum_{K\in \Th}\frac{1}{h} \int_K T_K |\lambda(\qq_K)| \ddd{x} \ddd{t} \leq \text{const.}
        \end{split}
        \end{equation}
    The last inequality holds due to the uniform bound of the density and a priori estimates derived in Section~\ref{sec:stability}.
     {Following the definition of $\lambda$ from \eqref{eq:lambda}, we have the bound $\lambda(\qq_K) \leq \lambda \leq \overline{u} + \overline{c}$, where $\overline{u} = \frac{\overline{m}}{\underline{\rho}}$ and $\overline{c} = \gamma \sqrt{\overline{p}/\underline{\rho}}$.}
        Thus the weak BV bound \eqref{eq:weakBV} on the diffusive terms holds.
    \end{proof}

Due to the flux correction term in the HTC scheme \eqref{eq:HTC_scheme}, an extra term given by
    \begin{align}
        \label{eq:correction_consistency}
        h^{d-1} \int_0^T \sum_{K\in \mathcal{T}}\sum_{\sigma \in \pK} \alpha_\sigma \jump{\pp_h(t)}_\sigma (P_h\phibf(t))_K \ddd{t}
    \end{align}
needs to be considered.
    The correction terms belonging to one arbitrary but fixed interface $\sigma = LK$ are
    \begin{equation}
        \label{eq:consistency-alpha-terms-interface}
        \frac{1}{h} \int_0^T \left(\alpha_\sigma \jump{\pp_h(t)}_\sigma \int_K \phibf(t) \ddd{x} - \alpha_\sigma \jump{\pp_h(t)}_\sigma \int_L \phibf(t) \ddd{x}\right) \ddd{t}.
    \end{equation}
    Analogously to the estimate for the diffusion terms in \cite[Section 9.3]{FeiLukMiz2019}, we consider an arbitrary but fixed point $\tilde x \in \sigma$. W.l.o.g. let $\tilde{x} = (\tilde{x}_s,x^\prime)$, where $x^\prime \in \mathbb{R}^{d-1}$, $s =1, \dots, d$.
    The Taylor series expansion for $x = (x_s,x^\prime) \in K$ and $x\in L$ with respect to the interface point $(\tilde{x}_s,x^\prime)$ yields, respectively,
    \begin{equation}
        \label{eq:Taylor}
        \phibf(x_s,x^\prime) = \phibf(\tilde{x}_s,x^\prime) \mp \xi \partial_s \phibf(\tilde{x}_s,x^\prime) + \mathcal{O}(h^2), \text{ where } \xi \in (0,h).
    \end{equation}
    Substituting the Taylor series expansion into \eqref{eq:consistency-alpha-terms-interface} we see that the terms multiplied with the leading order in the expansion $\phibf(\tilde{x}_s,x^\prime)$ vanish and the remaining terms can be estimated as follows
    \begin{align}
        \label{eq:AuxAlpha}
        \begin{split}
            &\left|\int_0^T \left(- \frac{1}{h} \alpha_\sigma\jump{\pp_h}_\sigma \int_0^h \int_\sigma \xi\partial_s \phibf(\tilde{x}_s,x^\prime)\ddd{\xi}\ddd{S_{x^\prime}} + \frac{1}{h} \alpha_\sigma\jump{\pp_h}_\sigma \int_0^h \int_\sigma -\xi\partial_s \phibf(\tilde{x}_s,x^\prime)\ddd{\xi}\ddd{S_{x^\prime}}\right)\right| \\
            &\leq \frac{2}{h} \int_0^T \left| \alpha_\sigma \jump{\pp_h}_\sigma \int_0^h \int_\sigma \xi\partial_s \phibf(\tilde{x}_s,x^\prime)\ddd{\xi}\ddd{S_{x^\prime}}\right| \ddd{t} \\
            & \lesssim h^d \int_0^T |\alpha_\sigma| \left|\jump{\pp_h}_\sigma\right|\ddd{t} ~\|\phibf\|_{C^1([0,T]\times\Omega)}.
        \end{split}
    \end{align}
Provided that we can prove a similar result as the weak BV property for $\jump{\qq_h}$ in Lemma \ref{lem:weakBV} on $\jump{\pp_h}$, the last term in \eqref{eq:AuxAlpha} tends to $0$ as $h \to 0^+$.
This property on $\jump{\pp_h}$ is subject to the next lemma.
\begin{lemma}
    { Let $\qq_h(t)$ be a solution of the HTC scheme \eqref{eq:HTC_scheme}.}
    { Then under Assumption \ref{ass:uniform_positivity} it holds}
        \begin{equation}
        \label{eq:BV_alpha}
        \int_0^T \sum_{\sigma \in \mathcal{E}} |\alpha_\sigma| \Big|\jump{\pp_h(t)}_\sigma\Big| h^d \ddd{t} \to 0 \quad \text{as} \quad h \to 0^+.
    \end{equation}
\end{lemma}
\begin{proof}
    The proof follows the lines of the proof of Lemma \ref{lem:weakBV}.
    Analogously to estimate \eqref{eq:Estimate_Pi}, we obtain for the correction term
    \begin{equation}
        \label{eq:Estimate_Pi_in_p}
        \int_0^T \sum_{\sigma \in \mathcal{E}} h^d \frac{\lambda_\sigma}{2} \frac{1}{T_K} \jump{\qq_h(t)}_\sigma {\cdot} \jump{\pp_h(t)}_\sigma \ddd{t} \geq \frac{\tilde\eta}{2}\int_0^T \sum_{\sigma \in \mathcal{E}} h^d \lambda_\sigma \frac{1}{T_K}\jump{\pp_h(t)}_\sigma  {\cdot} \jump{\pp_h(t)}_\sigma \ddd{t} \geq 0,
    \end{equation}
    where we have used that the inverse $(\HH_h)_\sigma^{-1}$ is strictly positive definite, thus there exists a {uniform} positive constant $\tilde \eta$ such that $(\HH_h)_\sigma^{-1} \geq \tilde{\eta} \mathbb{I}$.
    {This can be seen by taking the minimum over all eigenvalues of ${(\HH_h)_\sigma}^{-1}$ which are uniformly bounded away from zero due to Lemma \ref{lem:3}.}
    Consequently, $\jump{\qq_h}_\sigma \geq \tilde{\eta}\jump{\pp_h}_\sigma$.
    Note that the first term in \eqref{eq:Estimate_Pi_in_p} tends to zero due to \eqref{eq:h_limit_Pi}.
    Analogously to \eqref{eq:Aux3}, one obtains
    \begin{equation}
        \label{eq:Aux4}
        \int_0^T \sum_{\sigma \in \mathcal{E}} |\alpha_\sigma| \Big|\jump{\pp_h(t)}_\sigma\Big| h^d \ddd{t}
        \leq \left(\int_0^T  \sum_{\sigma \in \mathcal{E}} h^d T_K\frac{|\alpha_\sigma|^2}{\lambda_\sigma}\ddd{t}\right)^{1/2} \left(\int_0^T\sum_{\sigma \in \mathcal{E}}h^d\frac{\lambda_\sigma}{T_K} \Big|\jump{\pp_h(t)}\Big|^2\ddd{t}\right)^{1/2}.
    \end{equation}
    Then due to \eqref{eq:h_limit_Pi} and \eqref{eq:Estimate_Pi_in_p}, the second term on the right hand side in \eqref{eq:Aux4} tends to zero.
    The boundedness of the first term can be shown using the discrete trace inequality
    \begin{equation}
        \begin{split}
        \int_0^T  \sum_{\sigma \in \mathcal{E}} h^d T_K\frac{|\alpha_\sigma|^2}{\lambda_\sigma}\ddd{t} &\leq h \int_0^T \sum_{K\in \Th}\sum_{\sigma \in \partial K} \int_\sigma T_K\frac{|\alpha_\sigma|^2}{\lambda_\sigma} \ddd{S} \ddd{t} ~\\
        &{\leq h \int_0^T \sum_{K\in \Th} \frac{1}{h}\int_KT_K \frac{|\alpha|^2}{\tilde{\lambda}} \ddd{x}\ddd{t}}\leq \text{const.}
        \end{split}
    \end{equation}
    The last inequality holds due the uniform bounds on the density,  the well-posedness of $\alpha_\sigma$ and a priori estimates derived in Section \ref{sec:stability}.
    {In particular, the local diffusion coefficient $\lambda_\sigma$ can be estimated from below by $\tilde{\lambda}=\min_{K\in\Th} \min_k(|\lambda_k(\qq_L)|, |\lambda_k(\qq_K)|)$ at an arbitrary interface $\sigma = LK$. 
    Thus, it is uniformly bounded from below by $\underline{c} = \gamma \underline{p}/\overline{\rho}$ and 
    consequently $1/\tilde \lambda$ is bounded.
    With a similar argumentation, we obtain a uniform upper bound for $|\alpha| = \max_{K \in \Th}\alpha_\sigma(\qq_K,\qq_L)$ applying Lemma \ref{lem:Abgrall_alpha} in case of vanishing jumps in $\pp$. }
    This concludes the proof.
    \end{proof}

    It remains  to derive the estimate of the entropy production term in the entropy balance law \eqref{disc_entropy} which in the HTC formalism is discretized directly in the numerical scheme.
    We have {the following expression for the entropy production term}
    \begin{equation}
    h^{d-1}\int_0^T\sum_{K\in \Th}\frac{1}{2 T_K(t)}\sum_{\sigma \in \mathcal{E}}\lambda_\sigma\jump{\qq_h(t)}_\sigma {\cdot~} (\HH_h(t))_\sigma\jump{\qq_h(t)}_\sigma ~(P_h\varphi(t))_K\ddd{t}.
    \end{equation}
    The terms belonging to an arbitrary but fixed interface $\sigma = LK$ are given by
    \begin{equation}
        \label{eq:Aux2}
    \begin{split}
    \frac{1}{2h}\int_0^T &\left(\frac{1}{ T_L(t)}\lambda_\sigma\jump{\qq_h(t)}_\sigma  {\cdot~}  (\HH_h(t))_\sigma\jump{\qq_h(t)}_\sigma \int_L\varphi(t) \ddd{x}+\frac{1}{ T_K(t)}\lambda_\sigma\jump{\qq_h(t)}_\sigma  {\cdot~} (\HH_h(t))_\sigma\jump{\qq_h(t)}_\sigma \int_K\varphi(t) \ddd{x} \right)~\ddd{t}.
    \end{split}
    \end{equation}
    Similar as in the BV estimates for the local diffusion terms \eqref{eq:diffusion_consistency} and the consistency estimate for the correction terms \eqref{eq:correction_consistency}, we consider the Taylor series expansion \eqref{eq:Taylor} for $x = (x_s,x^\prime) \in K$ and $x\in L$, respectively.
    Substituting the Taylor series expansion into \eqref{eq:Aux2} and assuming w.l.o.g. that  $T_K \leq T_L$, we obtain the following estimate of \eqref{eq:Aux2}
    \begin{align*}
        \Big| \int_0^T&\left( \frac{1}{2h}\big((\Pi_h)_K+(\Pi_h)_L\big)\int_0^h\int_\sigma \varphi(\tilde{x}_s,x^\prime)\ddd{\xi}\ddd{S}_{x^\prime}+ \frac{1}{2h}\big(-(\Pi_h)_K+(\Pi_h)_L\big)\int_0^h\int_\sigma \xi\partial_s\varphi(\tilde{x}_s,x^\prime)\ddd{\xi}\ddd{S}_{x^\prime} \right)\Big|\\
        \leq& \frac{1}{h}\int_0^T (\Pi_h)_K \Big|\int_0^h \int_\sigma \varphi(\tilde{x}_s,x^\prime)\ddd{S}_{x^\prime}\Big| \ddd{t}+ \frac{1}{h}\int_0^T (\Pi_h)_K\Big|\int_0^h\int_\sigma \xi\partial_s\varphi(\tilde{x}_s,x^\prime)\ddd{\xi}\ddd{S}_{x^\prime} \Big|\ddd{t}\\
        \leq& h^d \int_0^T (\Pi_h)_K \ddd{t}~ \|\varphi\|_{C^1([0,T]\times\Omega)} + \frac{h^d}{2} \int_0^T (\Pi_h)_K \ddd{t}~ \|\varphi\|_{C^1([0,T]\times\Omega)} \to 0^+ \quad \text{for} \quad h \to 0^+,
    \end{align*}
    where we have used that $(\Pi_h)_K \geq 0$ for all $K \in \Th$ and estimate \eqref{eq:h_limit_Pi} to obtain the last convergence.

    \subsection{Consistency Formulation}
    Summarizing, we have shown the consistency of the HTC scheme \eqref{eq:HTC_scheme} for the Euler equations \eqref{eq:Euler_in_S}.

    \begin{theorem}[Consistency Formulation]\label{theo: consistency}
        Let $\qq_h^0$ be the initial data and $\qq_h$ a solution of the HTC scheme \eqref{eq:HTC_scheme} on the interval $[0,T]$ starting from $\qq_h^0$.

        Then, under Assumption \ref{ass:uniform_positivity} we have the following results for all $\tau\in(0,T)$:
        \begin{itemize}
            \item The approximate density $\rho_h$ fulfills
            \begin{equation} \label{98}
            \left[\int_\Omega \rho_h \varphi\ddd{x}\right]_{t = 0}^{t = \tau} = \int_0^\tau\int_\Omega \left( \rho_h \partial_t \varphi + \mm_h \cdot \nabla \varphi \right) \ddd{x}\ddd{t} + {\mathcal{O}}(h) \quad \text{as}\quad h \to 0
            \end{equation}
            for all $\varphi \in C^1([0,T], C^2(\Omega))$.
            \item The approximate momentum $\m_h$ fulfills
            \begin{equation} \label{99}
            \left[\int_\Omega \m_h \phibf\ddd{x}\right]_{t = 0}^{t = \tau} = \int_0^\tau\int_\Omega \left( \m_h \partial_t \phibf + \left(\frac{\m_h \otimes \m_h}{\rho_h}+ p_h \mathbb{I}\right) : \nabla \phibf \right) \ddd{x}\ddd{t} ~+~  {\mathcal{O}}(h)
            \end{equation}
            as $ h \to 0$ for all $\phibf \in C^1([0,T],C^2(\Omega; \mathbb{R}^d))$.
            \item The approximate entropy $S_h$ fulfills
            \begin{equation} \label{100}
            \left[\int_\Omega S_h \varphi\ddd{x}\right]_{t = 0}^{t = \tau} \geq \int_0^\tau\int_\Omega \left( S_h \partial_t \varphi + \frac{\mm_hS_h}{\rho_h} \cdot \nabla \varphi \right) \ddd{x}\ddd{t} +  {\mathcal{O}}(h) \quad \text{as}\quad h \to 0
            \end{equation}
        for all $\varphi \in C^1([0,T], C^2(\Omega))$.
        \item The total energy conservation holds
        \begin{equation} \label{101}
            \int_\Omega(\rho E)_h (\tau) \ddd{x} = \int_\Omega(\rho E)^0_h \ddd{x},
        \end{equation}
        where $ {\mathcal{O}}(h) \to 0$ as $h\to 0$.
        \end{itemize}
    \end{theorem}

\section{Convergence analysis via dissipative  weak solutions}\label{sec:convergence}
In the following we will establish the convergence of the HTC scheme \eqref{eq:HTC_scheme}.
As shown in the previous section, $\qq_h$ is a stable and consistent approximation of the complete Euler system.
We start by defining the notion of dissipative weak solutions developed in \cite{Feireisl2016} and the Young measures \cite{Ball} upon which we will prove the weak convergence to a dissipative weak solution and the strong convergence to a strong solution on the lifespan of the latter.

\subsection{Dissipative weak solutions of the Euler equations}

For completeness we repeat the definition of dissipative weak solutions of the Euler equations \eqref{eq:Euler_in_S} from \cite[Section 5.2.1]{Feireisl2021book}.

\begin{definition}[Dissipative weak solutions]\label{def:DMV}
Let the initial data satisfy
	\begin{eqnarray*}
&& \rho^0 \in L^\gamma(\Td), \ \m^0 \in L^\frac{2 \gamma}{\gamma +1}(\Td; \RR^d),\ \  S^0 \in L^\gamma(\Td)\\
&& E^0 =  E(\rho^0, \m^0, S^0)\, \mbox{ and }
\int_\Omega \rho^0 E^0  \, \dx < \infty.
	\end{eqnarray*}
	We say that $(\rho, \m, S)$ is a \emph{dissipative weak solution} to the Euler system
in $[0,T) \times \Td$,	$0 < T \leq \infty$, if the following holds:
	\begin{itemize}
		\item {\bf Regularity.}
        The solution $(\rho, \m, S)$ belongs to the class \footnote{\raggedright We use the standard notation, e.g., \newline $ C_{\rm weak, loc}([0,T); L^{{\nu}}(\Td)):=\{ {\nu} \in C(K; L^{{\nu}}(\Td)) \mbox{ for any compact set } K \subset [0,T) \mbox{ and } \int_{\Td} r(t,x) \varphi(x) \mbox{d} x \in C([0,T))\mbox{ for any } \varphi \in L^{{\nu}'}(\Td) \}.$  Here ${\nu}' > 1,$ $\frac{1} {{\nu}} + \frac{1}{{\nu}'}=1.$ Analogously, for $BV_{\rm weak}$ space we have $\int_{\Td} r(t,x) \varphi(x) \mbox{d} x \in BV([0,T)).$ }
		\begin{align*}
			\rho &\in C_{\rm weak, loc}([0,T); L^\gamma(\Td)),\
            \m \in C_{\rm weak, loc}([0,T); L^{\frac{2\gamma}{\gamma + 1}}(\Td; \RR^d)),\\
            S &\in L^\infty(0,T; L^\gamma(\Td)) \cap BV_{\rm weak}([0,T); L^\gamma(\Td)) \\\
			&\int_\Omega \rho E(\rho, \m, S)(t, \cdot) \, \dx \leq \int_\Omega \rho^0 E(\rho^0, \m^0, S^0) \, \dx \ \mbox{ for any }\ 0 \leq t < T.
			\nonumber
		\end{align*}
		\item {\bf Equation of continuity.}
		The integral identity
		\[
		\int_0^T \intTd{ \Big[ \rho \partial_t \varphi + \m \cdot \Grad \varphi \Big] } \dt =  - \intTd{ \rho^0 \varphi(0, \cdot)}
		\]
		holds for any $\varphi \in C^1_c([0,T) \times \Td)$.

		\item {\bf Momentum equation.}  The integral identity
		\begin{align}
			\int_0^T &\intTd{ \left[ \m \cdot \partial_t \bfphi + \mathds{1}_{\rho > 0} \frac{\m \otimes \m}{\rho} : \Grad \bfphi +
				p\left(\rho, \frac{S}{\rho}\right) \Div \bfphi \right] } \dt \nonumber  \\ &=
			- \int_0^T \int_{\overline{\Omega}} \Grad \bfphi : {\rm d}\mathfrak{R} (t) {\rm d}t  - \intTd{ \m^0 \cdot \bfphi(0, \cdot) }
			\label{E19}
		\end{align}
		holds for any $\bfphi \in C^1_c([0,T) \times \Td; \RR^d)$, where the Reynolds stress defect
		\begin{equation} \label{E20}
            \mathfrak{R} \in  L^\infty(0, T; \mathcal{M}^+(\overline{\Omega}; \RR^{d \times d}_{\rm sym})).\footnote{\raggedright
The space $\mathcal{M}^+(\overline{\Omega}; \RR^{d \times d}_{\rm sym})$ is the space of Radon measures ranging in the set of symmetric positive semi-definite matrices, i.e., $\mathcal{M}^+(\overline{\Omega}; \RR^{d \times d}_{\rm sym}):= \left\{\mu\in\mathcal{M}(\overline{\Omega}; \RR^{d \times d}_{\rm sym}),\,\int_{\overline{\Omega}}\phi(\xi\otimes\xi):\mbox{d}\mu \geq 0 \;\text{for all} \; \xi\in\Bbb R^d, \; \phi\in C(\overline{\Omega}), \; \phi\geq 0\right\}$.}
		\end{equation}

\item{\bf Entropy inequality.}
\begin{align} \label{DS10}
& \left[ \intTd{ S \varphi } \right]_{t = \tau_1-}^{t = \tau_2+} \geq
\int_{\tau_1}^{\tau_2} \intTd{ \left[ S \partial_t \varphi + \left< {\mathcal{V}_{t,x}}; 1_{\tvr > 0} \left( \tvS \frac{\tvm}{\tvr} \right) \right> \cdot \Grad \varphi \right] } \dt, \\
& \hspace{8cm} S(0-, \cdot) = S^0, \nonumber
\end{align}
for any $0 \leq \tau_1 \leq \tau_2 < T$, any $\varphi \in C^1_c([0,T) \times \Td)$, $\varphi \geq 0$,
where
$
\{ {\mathcal{V}}_{t,x} \}_{ (t,x) \in (0,T) \times \Td }
$
is a parametrized probability measure (Young measure),
\[
{\mathcal{V}_{t,x}} \in L^\infty((0,T) \times \Td); {\mathcal{P}}(\Bbb R^{d + 2})),\ \Bbb R^{d+2} = \left\{ \tvr \in \Bbb R, \tvm \in \Bbb R^d, \tvS \in \Bbb R \right\};
\]
\begin{equation} \label{DS10a}
\left< {\mathcal{V}}; \tvr \right> = \rho, \ \left< {\mathcal{V}}; \tvm \right> = \m,\
\left< {\mathcal{V}}; \tvS \right> = S.
\end{equation}

		\item{\bf Compatibility of the energy and Reynolds stress defects.} There exists a non--increasing function $\mathcal{E}:
		[0,T) \to [0, \infty)$ satisfying
		\begin{align}
			\mathcal{E}(0-) &= \intTd{ \rho^0 E(\rho^0, \m^0, S^0) }, \nonumber \\
			\mathcal{E}(\tau+) & =  \intTd{ \rho E(\rho, \m, S)(\tau, \cdot) } + \mathfrak{E} 	\quad \mbox{ for any } 0 \leq \tau < T,
			\label{E21}
		\end{align}
		where $\mathfrak{E} \in   L^\infty(0, T; \mathcal{M}^+(\Td)) $ is the energy defect satisfying
		\begin{align}
		 \min \left\{ 2, d(\gamma - 1) \right\} \mathfrak{E}  \leq  {\rm trace} [\mathfrak{R}] \leq  \max \left\{ 2, d(\gamma - 1) \right\}  \mathfrak{E}.
			\label{eq:def}
		\end{align}
		\end{itemize}

\end{definition}

Let  $\left\lbrace\qq_h\right\rbrace_{h>0}:=\left\lbrace\rho_h, \m_h, S_h\right\rbrace_{h>0}$ be a family of numerical solutions computed with the HTC scheme \eqref{eq:HTC_scheme}.
Note that due to Theorem \ref{theo: consistency} $\left\lbrace\rho_h, \m_h, S_h\right\rbrace_{h>0}$ is a consistent approximation of the Euler
equations~\eqref{eq:Euler_in_S}.
Following \cite[Section 9.24]{Feireisl2021book}, $\left\lbrace\qq_h\right\rbrace_{h>0}$  generates a Young  measure $\mathcal{V}_{t,x}$
and  converges up to a subsequence to a dissipative weak solution $(\rho, \m, S)$  of the Euler equations in the sense of Definition \ref{def:DMV}.

\begin{theorem}[Weak convergence]\label{theo:DMVsolutions}
    Let the initial data $\qq^0$ satisfy \eqref{eq:init_data} and \eqref{eq:init_data_thermodyn}.
    Let further $\left\lbrace\qq_h\right\rbrace_{h>0} =\left\lbrace(\rho_h, \m_h, S_h)\right\rbrace_{h>0}$ be a family of numerical solutions generated by the HTC scheme \eqref{eq:HTC_scheme}.
    Further suppose that Assumption \ref{ass:uniform_positivity} holds, i.e. that there exist constants $\underline{\rho} > 0$ and $\overline{E} > 0$ such that
    \begin{equation}
        0 < \underline{\rho} \leq \rho_h, \quad (\rho E)_h \leq \overline{E} \quad \text{uniformly for } \quad h \to 0.
    \end{equation}

Then there exists a subsequence $\left\lbrace\qq_{h_n}\right\rbrace_{{h_n}>0}$, such that
    \begin{equation}
       (\rho_{h_n},\mm_{h_n},S_{h_n}) \to (\rho, \m, S) \equiv \left\langle\DMV;(\widetilde{\rho},\widetilde{\m},\widetilde{S})\right\rangle
    \end{equation}
    weakly-$\ast$ in $L^\infty((0,T)\times\Omega;\mathbb{R}^{d+2}).$

Furthermore, the defects are
    \begin{align}
    \label{eq:defects}
        &\mathfrak{C} =  \left\langle\DMV; \widetilde\rho E (\widetilde{\rho},\widetilde{\m},\widetilde{S})\right\rangle - \rho E(\rho,\mm,S), \\
        &\mathfrak{R} =  \left\langle\DMV;\frac{\widetilde\mm\otimes\widetilde\mm}{\widetilde\rho}+p\left(\widetilde\rho,\frac{\widetilde{S}}{\widetilde\rho}\right)\mathbb{I}\right\rangle
    - \frac{\mm\otimes\mm}{\rho}-p\left(\rho,\frac{S}{\rho}\right)\mathbb{I},
    \end{align}
where
\begin{align*}
        &\rho_{h_n}  E(\rho_{h_n},\mm_{h_n},S_{h_n}) \to  \left\langle\DMV; \widetilde\rho E (\widetilde{\rho},\widetilde{\m},\widetilde{S})\right\rangle \text{weakly}-\ast \text{ in } L^\infty( (0,T) \times \Omega), \ \mbox{ as } h_n \to 0,\\
    &\frac{\mm_{h_n}\otimes\mm_{h_n}}{\rho_{h_n}}+p\left(\rho_{h_n},\frac{S_{h_n}}{\rho_{h_n}} \right)\mathbb{I} \to
        \left\langle\DMV;\frac{\widetilde\mm\otimes\widetilde\mm}{\widetilde\rho}+p\left(\widetilde\rho,\frac{\widetilde{S}}{\widetilde\rho}\right)\mathbb{I}\right\rangle \\
        &\text{weakly}-\ast \text{ in } L^\infty( (0,T) \times \Omega; \Bbb R^{d \times d})) \  \mbox{ as } h_n \to 0.
\end{align*}

\end{theorem}
\begin{proof}
    Under Assumption \ref{ass:uniform_positivity}, Lemma \ref{lem:3} gives the $L^\infty$ bounds \eqref{eq:apriori_estimates} and
    \begin{equation*}
        \frac{\m_h\otimes\m_h}{\rho_h} \in L^\infty((0,T)\times \Omega), \quad S_h \frac{\m_h}{\rho_h}  \in L^\infty((0,T)\times\Omega).
    \end{equation*}
    Applying the Fundamental Theorem on Young Measures, see \cite{Ball} or \cite{Pedregal1997}, it implies the existence of a convergent subsequence and a parametrized probability measure $\left\lbrace\DMV\right\rbrace_{(t,x)\in (0,T)\times\Omega}$ such that $(\rho_{h_n},\mm_{h_n},S_{h_n}) \to (\rho, \m, S) \equiv \left(\left\langle\DMV;\widetilde{\rho}\right\rangle,\left\langle\DMV;\widetilde{\mm}\right\rangle,\left\langle\DMV;\widetilde{S}\right\rangle\right) $ weakly-$\ast$ in $L^\infty((0,T)\times\Omega)$ as $h_n \to 0$.
    Moreover, the following terms
    \begin{equation*}
        \frac{\m_{h_n}\otimes\m_{h_n}}{\rho_{h_n}}, \quad p_{h_n} = p\left(\rho_{h_n},\frac{S_{h_n}}{\rho_{h_n}}\right), \quad S_{h_n} \frac{\m_{h_n}}{\rho_{h_n}}, \quad \rho_{h_n} E\left(\rho_{h_n},\m_{h_n},S_{h_n}\right)
    \end{equation*}
converge   weakly-$\ast$   to
    \begin{equation*}
        \left\langle\DMV;\frac{\widetilde{\m}\otimes\widetilde{\m}}{\widetilde{\rho}}\right\rangle, \quad \left\langle\DMV;p\left(\widetilde{\rho},\frac{\widetilde{S}}{\widetilde{\rho}}\right)\right\rangle, \quad \left\langle\DMV;\widetilde{S} \frac{\widetilde{\m}}{\widetilde{\rho}}\right\rangle,\quad \left\langle\DMV;(\rho E) (\widetilde{\rho},\widetilde{\m},\widetilde{S})\right\rangle
    \end{equation*}
    in $L^\infty((0,T)\times\Omega)$, respectively.
    Consequently, passing to the limit as $h \to 0$ in the consistency formulation \eqref{99} - \eqref{101} we obtain
    \begin{equation*}
 	\int_0^T \intTd{ \Big[ \rho \partial_t \varphi + \m \cdot \Grad \varphi \Big] } \dt =  - \intTd{ \rho^0 \varphi(0, \cdot)}
    \end{equation*}
    for all $\varphi \in C^1_c([0,T), C^2(\Omega))$ for the density and
\begin{eqnarray}
			&& \int_0^T  \intTd{ \left[ \m \cdot \partial_t \bfphi + \mathds{1}_{\rho > 0} \frac{\m \otimes \m}{\rho} : \Grad \bfphi +
				p\left(\rho, \frac{S}{\rho}\right) \Div \bfphi \right] } \dt \nonumber  \\   \nonumber
              &&= - \int_0^T \intTd{ \Grad \bfphi : {\rm d}\mathfrak{R} (t) }  - \intTd{ \m^0 \cdot \bfphi(0, \cdot) }
\end{eqnarray}
      for all $\bfphi \in C^1_c([0,T); C^2(\Omega; \mathbb{R}^d))$ for the momentum.
    Further, we have
\begin{eqnarray*}
 && \left[ \intTd{ S \varphi } \right]_{t = \tau_1-}^{t = \tau_2+} \geq
\int_{\tau_1}^{\tau_2} \intTd{ \left[ S \partial_t \varphi + \left< {\mathcal{V}_{t,x}}; 1_{\tvr > 0} \left( \tvS \frac{\tvm}{\tvr} \right) \right> \cdot \Grad \varphi \right] } \dt, \\
&& \hspace{8cm} S(0-, \cdot) = S^0
\end{eqnarray*}
    for any $0 \leq \tau_1 \leq \tau_2 < T $ and all $\varphi \in C^1_c([0,T); C^2(\Omega))$ for the entropy
    and
    \begin{equation}
        \int_\Omega \rho(t, \cdot)  E(\rho, \m, S) (t, \cdot) \ddd{x}  + \int_{\overline{\Omega}} \ddd \mathfrak{E} (t) \leq  \int_\Omega \rho^0 E^0 \, \ddd{x} \quad \mbox{ for a.a. } t \in (0,T)
    \end{equation}
    for the total energy. The compatibility condition \eqref{eq:def} follows by comparing the energy and trace of the Reynolds defect.
    Applying the density argument on the test function, we finally get that ${(\rho, \m, S)}$ is a dissipative weak solution of the complete Euler system \eqref{eq:Euler_in_S} which concludes the proof.
\end{proof}

\subsection{Strong convergence}

In the above theorem it was proven that the numerical solutions of the HTC scheme converge weakly-$\ast$ to a dissipative weak solution.
Employing dissipative weak-strong uniqueness results from \cite{Brezina2018}, we can show the strong convergence to a strong solution as long as it  exists.

\begin{theorem}\label{theo:strongconv}
Let assumptions of Theorem~\ref{theo:DMVsolutions} hold.
In addition, let the Euler equations \eqref{eq:Euler_in_S} possess the unique strong (continuously differentiable) solution $\qq = (\rho, \m, S)$ emanating from the initial data \eqref{eq:init_data} and \eqref{eq:init_data_thermodyn}.

    Then
    \begin{equation}
        \qq_h \to \qq \quad \text{strongly} \quad \text{in} \quad L^q((0,T)\times\Omega{; \mathbb{R}^{d+2}}), \quad q \in [1,\infty).
    \end{equation}
    More precisely,
    \begin{align*}
        &\rho_h \to \rho &&\text{weakly-$\ast$ in} \quad L^\infty((0,T)\times\Omega) &&\text{and strongly in} \quad L^q((0,T)\times \Omega)\\
        &\m_h \to \m &&\text{weakly-$\ast$ in} \quad L^\infty((0,T)\times\Omega)  &&\text{and strongly in} \quad L^q((0,T)\times \Omega; \mathbb{R}^{d})\\
        &S_h \to S &&\text{weakly-$\ast$ in} \quad L^\infty((0,T)\times\Omega) &&\text{and strongly in} \quad L^q((0,T)\times \Omega) .
    \end{align*}
\end{theorem}

In numerical experiments, the justification and visualization of the weak convergence of a numerical scheme is not practical.
Thus as done in \cite{ALO2023,Lukacova2023}, we apply $\mathcal{K}-$ convergence introduced in \cite{Feireisl2021,Feireisl2021book}.
As shown in \cite[Theorem 10.5]{Feireisl2021book}, or \cite{Feireisl2021}, the Ces\'aro averages and first variances converge strongly to a dissipative weak solution even if the strong solution does not exist and there are infinitely many (dissipative)-weak solutions.
Applying the techniques  from \cite[Theorem 7.9]{Feireisl2021book} and \cite{Feireisl2021}, we obtain the same result for the HTC scheme. We have the strong convergence of the Ces\'aro averages
\begin{equation}
    \frac{1}{N}\sum_{n=1}^{N} \qq_{h_n} \to ({\rho, \m, S}) \quad \text{ in } \quad L^q((0,T)\times\Omega; ~\mathbb{R}^{d+2}),
\end{equation}
for $N\to\infty$ and any $q\in [1,\infty)$, and the $L^1$ convergence of the first variance
\begin{equation}
    \frac{1}{N}\sum_{n=1}^{N} \left|\qq_{h_n} - \frac{1}{N}\sum_{k=1}^{{N}} \qq_{h_k} \right|
    \to 0 \quad \text{in} \quad L^1((0,T)\times \Omega) \quad \text{for} \quad N \to \infty.
\end{equation}

\section{Numerical results}\label{sec:NumRes}
In this section we illustrate the theoretical results proven in the previous sections.
Note that the space semi-discrete version of the HTC scheme discussed in this work preserves the total energy in space.
Thus, the space semi-discrete HTC scheme \eqref{eq:HTC_scheme} is integrated in time with a fifth order explicit Runge Kutta method to reduce the time error on the total energy with a very high accuracy.
Note that in the fully implicit version, the total energy preservation in space-time can be guaranteed, see \cite{Abgrall2023HTC}.
The overall accuracy of the HTC scheme remains first order for smooth solutions, since the space discretization as given above is first order.
In the following we focus on the numerical assessment of the weak, strong and $\mc{K}$-convergence of the HTC method at the example of the Kelvin-Helmholtz (KH) problem.

We consider the set-up from \cite{Feireisl2021book} which was also used in the context of discontinuous Galerkin schemes \cite{ALO2023}.
The initial condition for the KH problem on the unit box $~\Omega = [0,1]\times [0,1]$ is given by
\begin{equation}
    \label{eq:KH}
    \left(\rho, v_1, v_2, p\right) =
    \begin{cases}
        \left(2, -0.5, 0, 2.5\right) &\text{if } I_1 < x_2 < I_2 \\
        \left(1,  0.5, 0, 2.5\right) &\text{else.}
    \end{cases}
\end{equation}
The interface profiles, denoted by $I_j = I_j(\x)$ for $j=1,2$, are defined as $I_j = J_j + \varepsilon Y_j(\x)$ which generate small perturbations around the constant horizontal lines $J_1 = 0.25$ and $J_2 = 0.75$.
Further, we set
\begin{equation}
    Y_j(\x) = \sum_{k=1}^m a_j^k \cos(b_j^k + 2 k \pi x_1), \quad j= 1,2,
\end{equation}
where the parameters $a_j^k \in [0,1]$ and $b_j^k$, $j=1,2$, $k=1, \dots, m$ are fixed.
The coefficients $a_j^k$ are normalized such that $\sum_{k=1}^m a_j^k = 1$ to guarantee that the perturbation around $J_j$ are of the order of $\varepsilon$, thus $|I_j(\x) - J_j| \leq \varepsilon$ for $j=1,2$.
In the numerical simulation we use $m=10$ modes and the perturbation $\varepsilon = 0.01$, as done in \cite{Feireisl2021book}.
Further we apply periodic boundary conditions in both directions.

The HTC scheme directly evolves the state variables $ \qq_h = (\rho_h, \mm_h, S_h)$ in which the convergence has been proven theoretically. Moreover, HTC numerical solution preserves discrete energy conservation of the total energy $(\rho E)_h.$
Note that the KH initial data as defined in \eqref{eq:KH} are discontinuous.

The experimental order of convergence (EOC) is computed based on the Ces\'aro averages and the first variance which are given by
\begin{equation}
    \tilde{U}_{h_n} = \frac{1}{n} \sum_{j=1}^n U_{h_j}, \quad U_{h_n}^\dagger = \frac{1}{n} \sum_{j=1}^n \Big|U_{h_j} - \tilde{U}_{h_n}\Big|,
\end{equation}
for $U \in \left\{\rho, \mm, S, \rho E\right\}$.
To compute the errors, we use the standard $L^1$ norm on the difference between the numerical and reference solution given by $E_1$.
Analogously, $L^1$ norm is applied to the Ces\'aro averages and the first variance denoted by $E_2$ and $E_3$, respectively.
Hence they are defined as
\begin{equation}
    E_{1_n} = \left\|U_{h_n} - U_{h_N}\right\|, \quad E_{2_n} = \left\|\tilde{U}_{h_n} - \tilde{U}_{h_N}\right\|, \quad E_{3_n} = \left\|U^\dagger_{h_n} - U^\dagger_{h_N}\right\|,
\end{equation}
where $h_n = \frac{1}{n}$, $n$ is the number of mesh cells in one direction and a finest mesh resolution has  $N \times N$ cells.
Based on the different errors $E_1, E_2, E_3$ defined above, the EOC is obtained via
\begin{equation}
    \text{EOC}_k = \log_2\left(E_{k_n}/E_{k_{2n}}\right), \quad k = 1, 2, 3,
\end{equation}
where subsequent meshes are refined by a factor of two.

\begin{table}[!htb]
    \centering
    \renewcommand{\arraystretch}{1.25}
    \begin{tabular}{cl|cr|cr|cr}
        &\multirow{2}{*}{$h_n$} & \multicolumn{2}{c|}{$E_{1_n}$}&\multicolumn{2}{c|}{$E_{2_n}$}&\multicolumn{2}{c}{$E_{3_n}$}\\
       & & error & EOC$_1$ & error & EOC$_2$ & error & EOC$_3$ \\\hline\hline
        \multirow{5}{*}{$\rho$}&  1/256 &  1.896E-1 &  ---   & 1.5297E-1 &  ---   & 1.351E-1 &   ---  \\
        &  1/512 &  1.609E-1 &  0.237 & 9.3815E-2 &  0.705 & 7.774E-2 &  0.798 \\
        &  1/1024 &  1.654E-1 & -0.040 & 6.1914E-2 &  0.600 & 5.008E-2 &  0.634 \\
        & 1/2048 &  2.062E-1 & -0.318 & 4.1086E-2 &  0.592 & 3.207E-2 &  0.643 \\
        & 1/4096 &  1.916E-1 &  0.106 & 2.4145E-2 &  0.767 & 1.716E-2 &  0.902 \\
\hline\hline
        \multirow{5}{*}{$\rho v_1$} & 1/256  &  1.912E-1 &  ---   & 1.328E-1 &  ---   & 1.100E-1 & ---   \\
        & 1/512  &  1.309E-1 &  0.546 & 6.783E-2 &  0.969 & 5.568E-2 &  0.983 \\
        & 1/1024  &  1.770E-1 & -0.434 & 4.873E-2 &  0.477 & 3.643E-2 &  0.612 \\
        & 1/2048 &  1.672E-1 &  0.082 & 3.328E-2 &  0.550 & 2.350E-2 &  0.632 \\
        & 1/4096 &  1.693E-1 & -0.018 & 2.179E-2 &  0.611 & 1.441E-2 &  0.705 \\
\hline\hline
        \multirow{5}{*}{$\rho v_2$} &  1/256 &  1.902E-1 &  ---   &  1.514E-1 &  ---   & 1.226E-1 &  ---   \\
        &  1/512 &  1.318E-1 &  0.529 &  8.711E-2 &  0.798 & 6.563E-2 &  0.902 \\
        &  1/1024 &  1.777E-1 & -0.431 &  4.775E-2 &  0.867 & 4.069E-2 &  0.689 \\
        & 1/2056 &  1.780E-1 & -0.002 &  3.209E-2 &  0.573 & 2.611E-2 &  0.640 \\
        & 1/4096 &  2.118E-1 & -0.251 &  2.280E-2 &  0.493 & 1.348E-2 &  0.953 \\
\hline\hline
        \multirow{5}{*}{$S$} &1/256  & 1.780E-1 &  ---   & 1.469E-1 &  ---   & 1.279E-1  & ---   \\
        &1/512  &1.472E-1 &  0.274 & 8.997E-2 &  0.708 & 7.673E-2  & 0.738 \\
        &1/1024  &1.458E-1 &  0.013 & 5.927E-2 &  0.602 & 4.836E-2  & 0.666 \\
        &1/2048 &1.852E-1 & -0.345 & 3.840E-2 &  0.626 & 3.126E-2  & 0.629 \\
        &1/4096 &1.755E-1 &  0.077 & 2.249E-2 &  0.772 & 1.683E-2  & 0.893 \\
\hline\hline
        \multirow{5}{*}{$\rho E$} &1/256  & 2.693E-1 &  ---   & 0.201E-1 &  ---   & 1.801E-1 &  ---   \\
        &1/512  & 2.307E-1 &  0.223 & 0.132E-1 &  0.602 & 1.006E-1 &  0.840  \\
        &1/1024  & 2.727E-1 & -0.241 & 0.859E-2 &  0.626 & 6.601E-2 &  0.608  \\
        &1/2048 & 2.784E-1 & -0.030 & 0.553E-2 &  0.636 & 4.603E-2 &  0.520  \\
        &1/4096 & 2.830E-1 & -0.024 & 0.331E-2 &  0.739 & 2.512E-2 &  0.874  \\
\hline\hline
    \end{tabular}
    \caption{Kelvin-Helmholtz problem: Convergence study using the errors $E_1, E_2, E_3$ with {$N=8192$} applied on the density, momentum, total entropy and total energy. }
    \label{tab:EOC}
\end{table}

\begin{figure}[h!]

    \vspace{-5mm}

    \begin{subfigure}{0.5\textwidth}
        \hspace{-3.5cm}\includegraphics[scale=0.18]{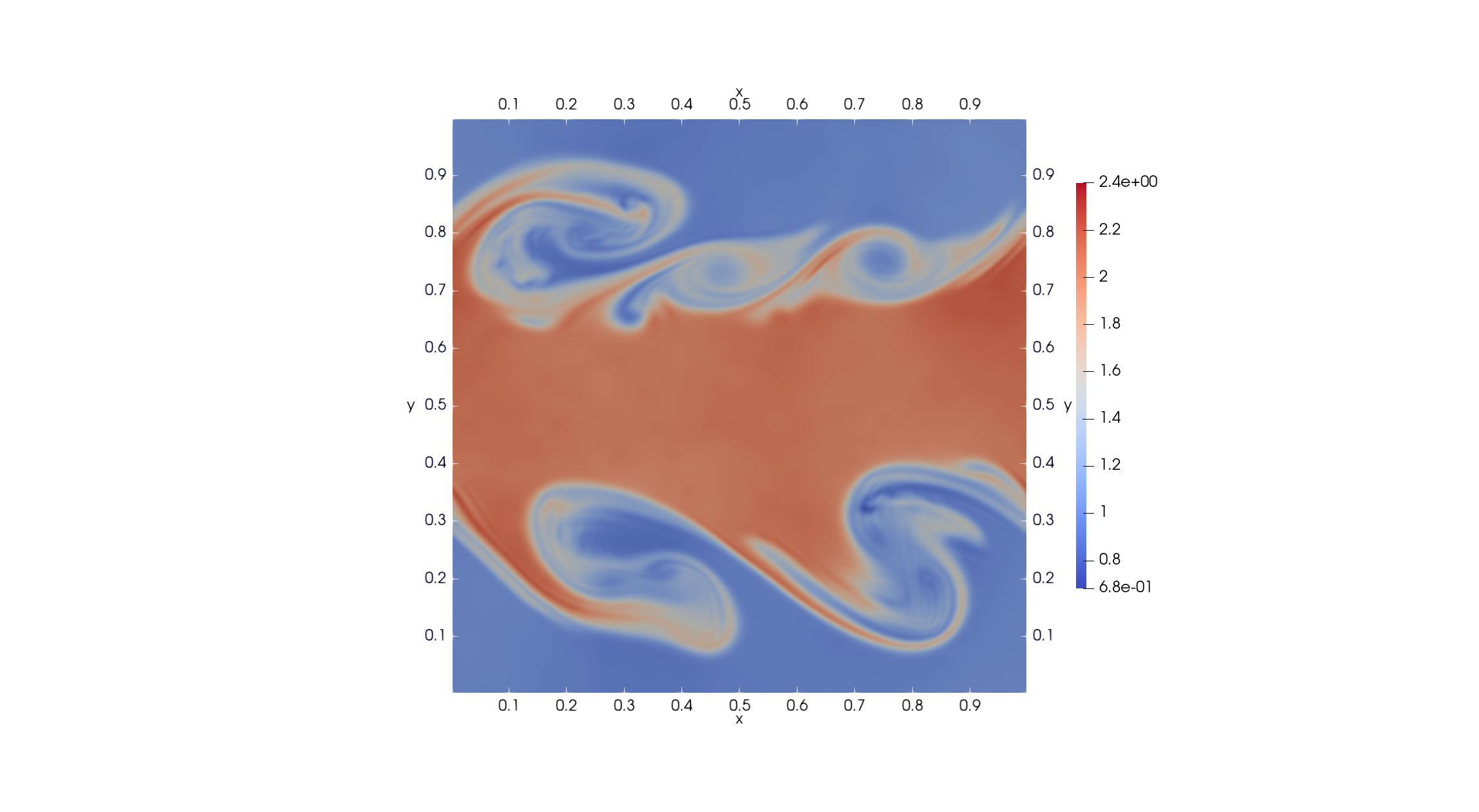}
        \vspace{-1.5cm}
        \caption{$n=256$}
        \label{fig:KH_256}
    \end{subfigure}
    \begin{subfigure}{0.5\textwidth}
       \hspace{-3.cm} \includegraphics[scale=0.18]{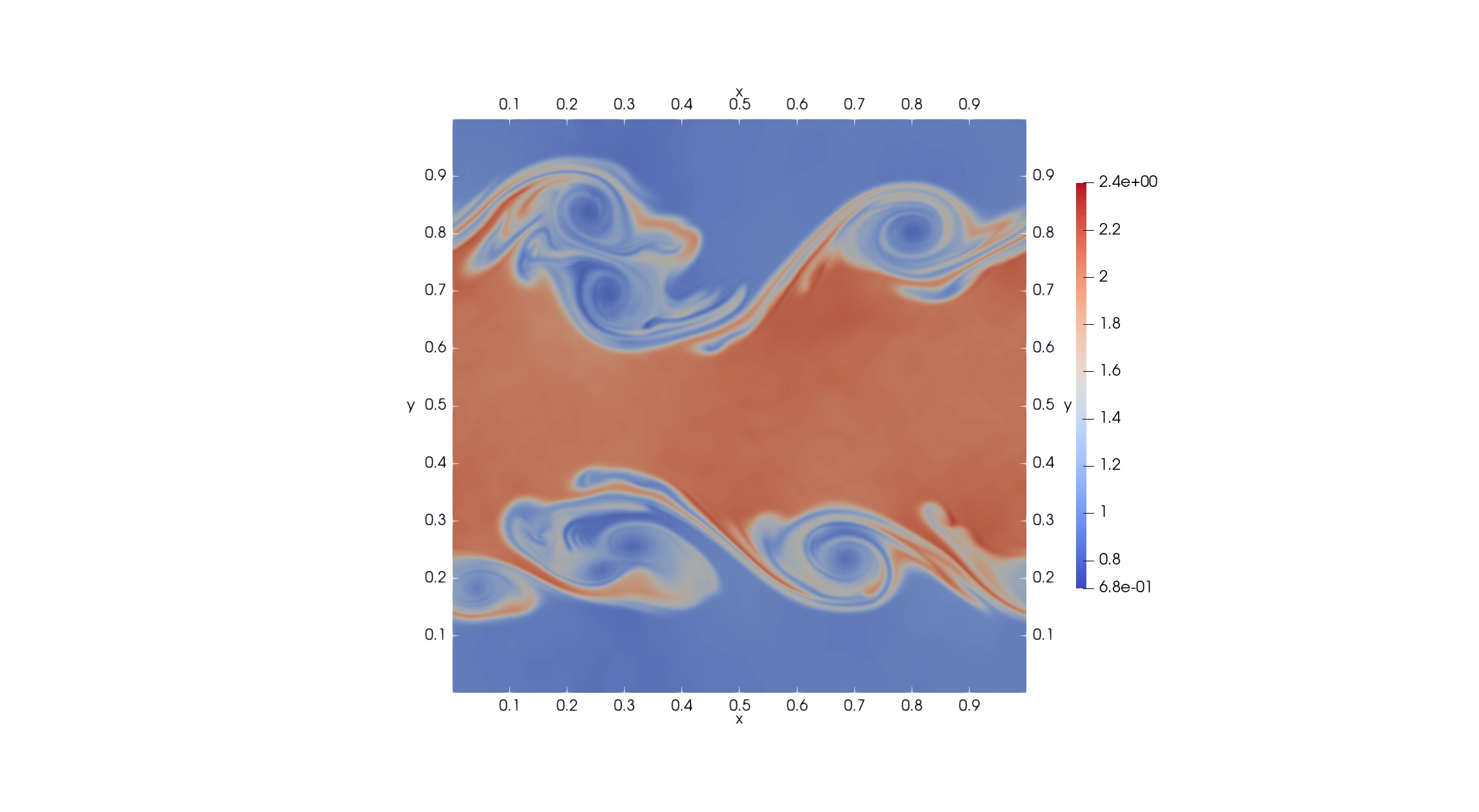}
       \vspace{-1.5cm}
        \caption{$n=512$}
        \label{fig:KH_512}
    \end{subfigure}

    \vspace{-0.7cm}

    \begin{subfigure}{0.5\textwidth}
         \hspace{-3.5cm} \includegraphics[scale=0.18]{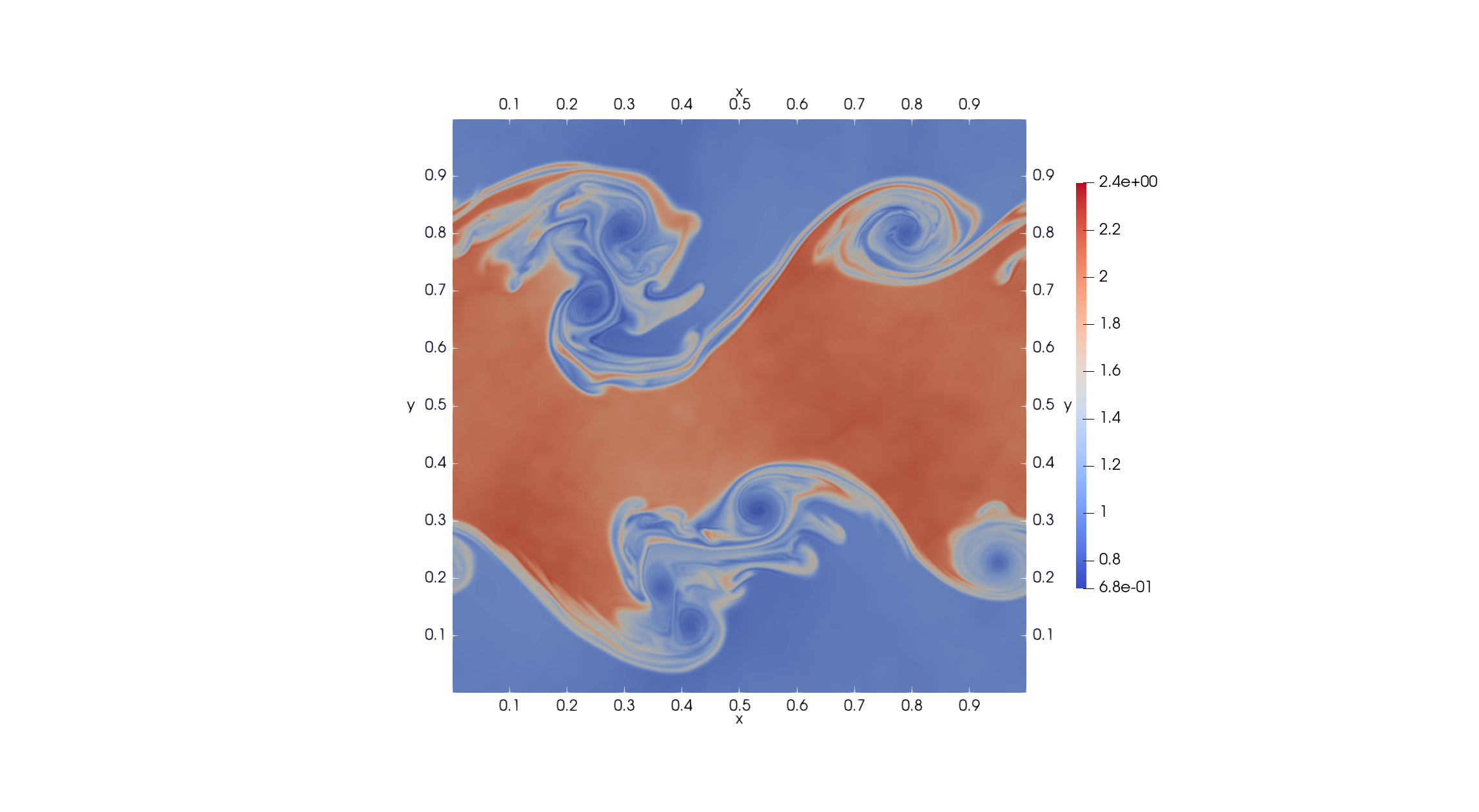}
         \vspace{-1.5cm}
        \caption{$n=1024$}
        \label{fig:KH_1024}
    \end{subfigure}
    \begin{subfigure}{0.5\textwidth}
        \hspace{-3.cm}  \includegraphics[scale=0.18]{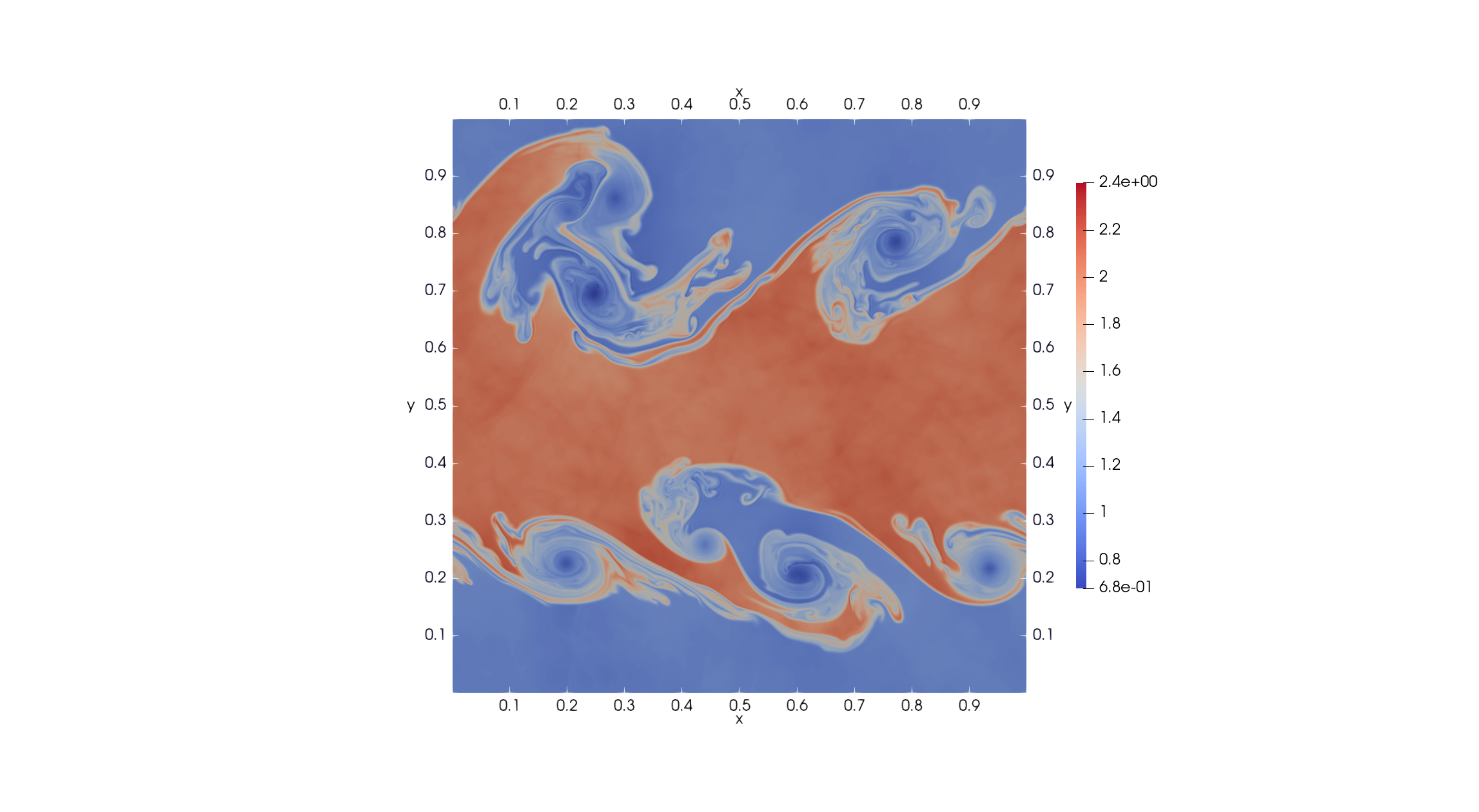}
        \vspace{-1.5cm}
        \caption{$n=2048$}
        \label{fig:KH_2048}
    \end{subfigure}

    \vspace{-0.7cm}

        \begin{subfigure}{0.5\textwidth}
        \hspace{-3.5cm } \includegraphics[scale=0.18]{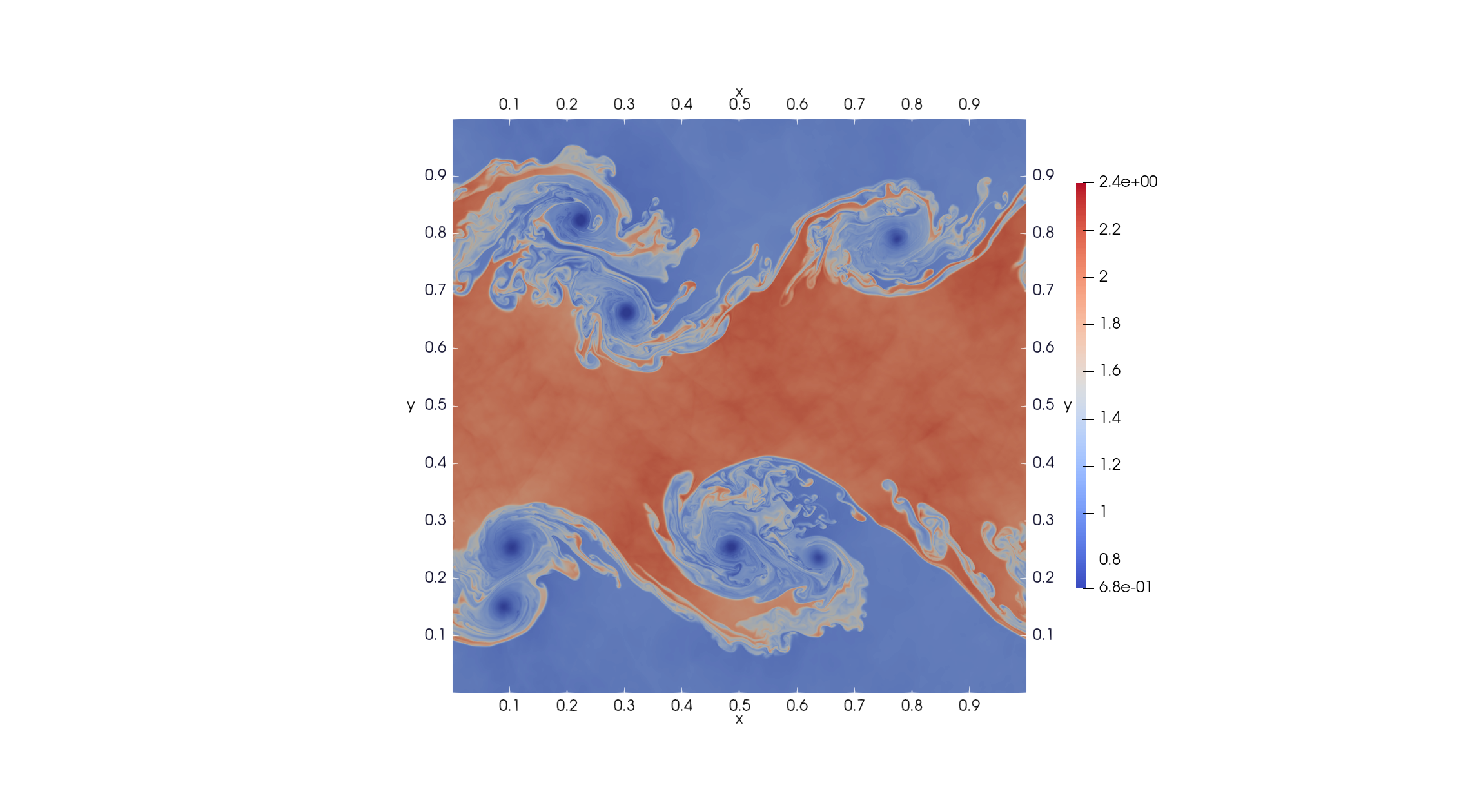}
        \vspace{-1.5cm}
        \caption{{$n=4096$}}
        \label{fig:KH_4086}
    \end{subfigure}
        \begin{subfigure}{0.5\textwidth}
       \hspace{-3.cm } \includegraphics[scale=0.18]{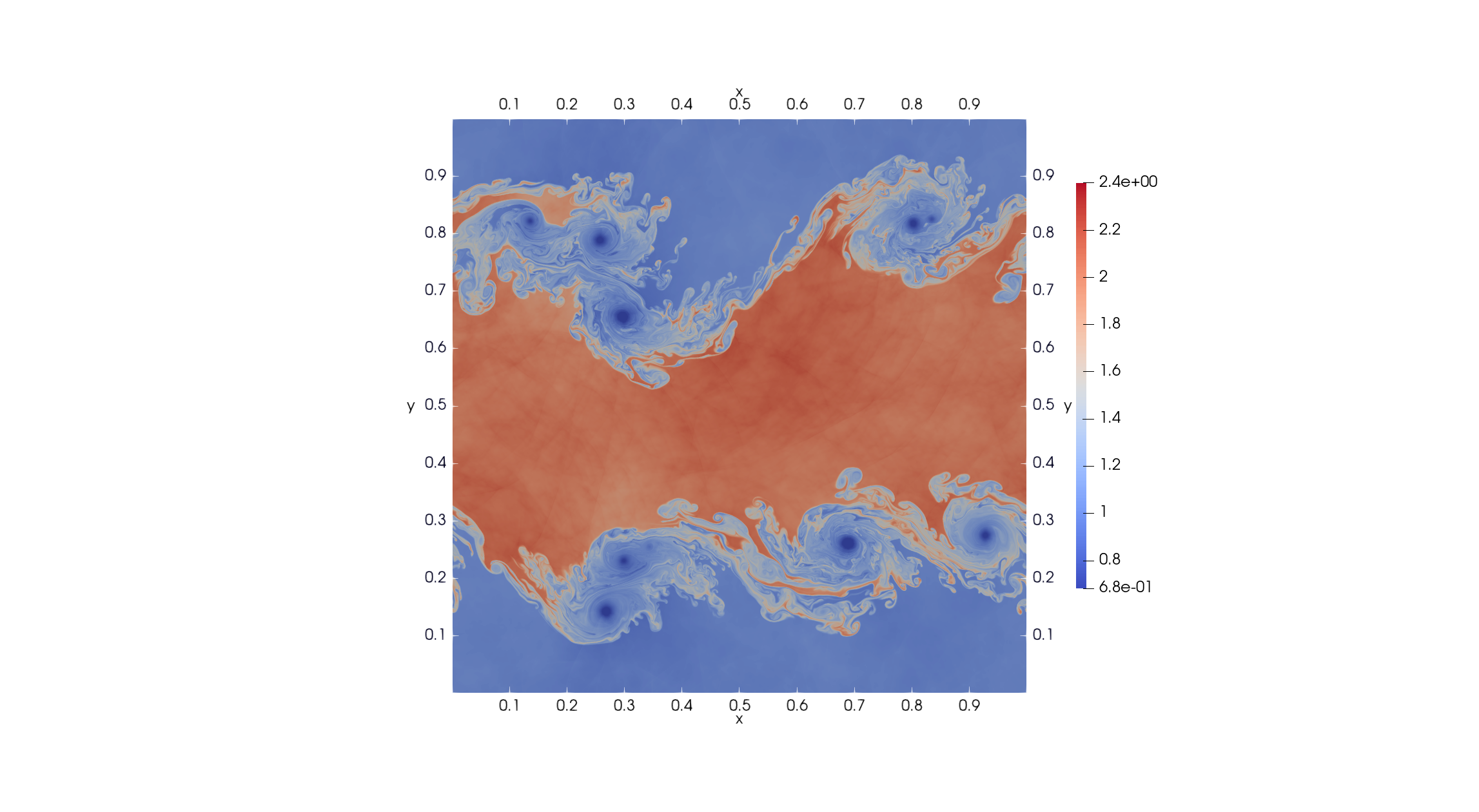}
       \vspace{-1.5cm}
        \caption{{$n=8192$}}
        \label{fig:KH_8192}
    \end{subfigure}
    \vspace{-0.2cm}
    \caption{Kelvin-Helmholtz problem: Density computed with the HTC scheme at time $T=2$ for different meshes with $n \times n$ cells.}
    \label{fig:KH}
\end{figure}

\begin{figure}[h!]

    \vspace{-5mm}

    \begin{subfigure}{0.5\textwidth}
        \includegraphics[scale=0.27]{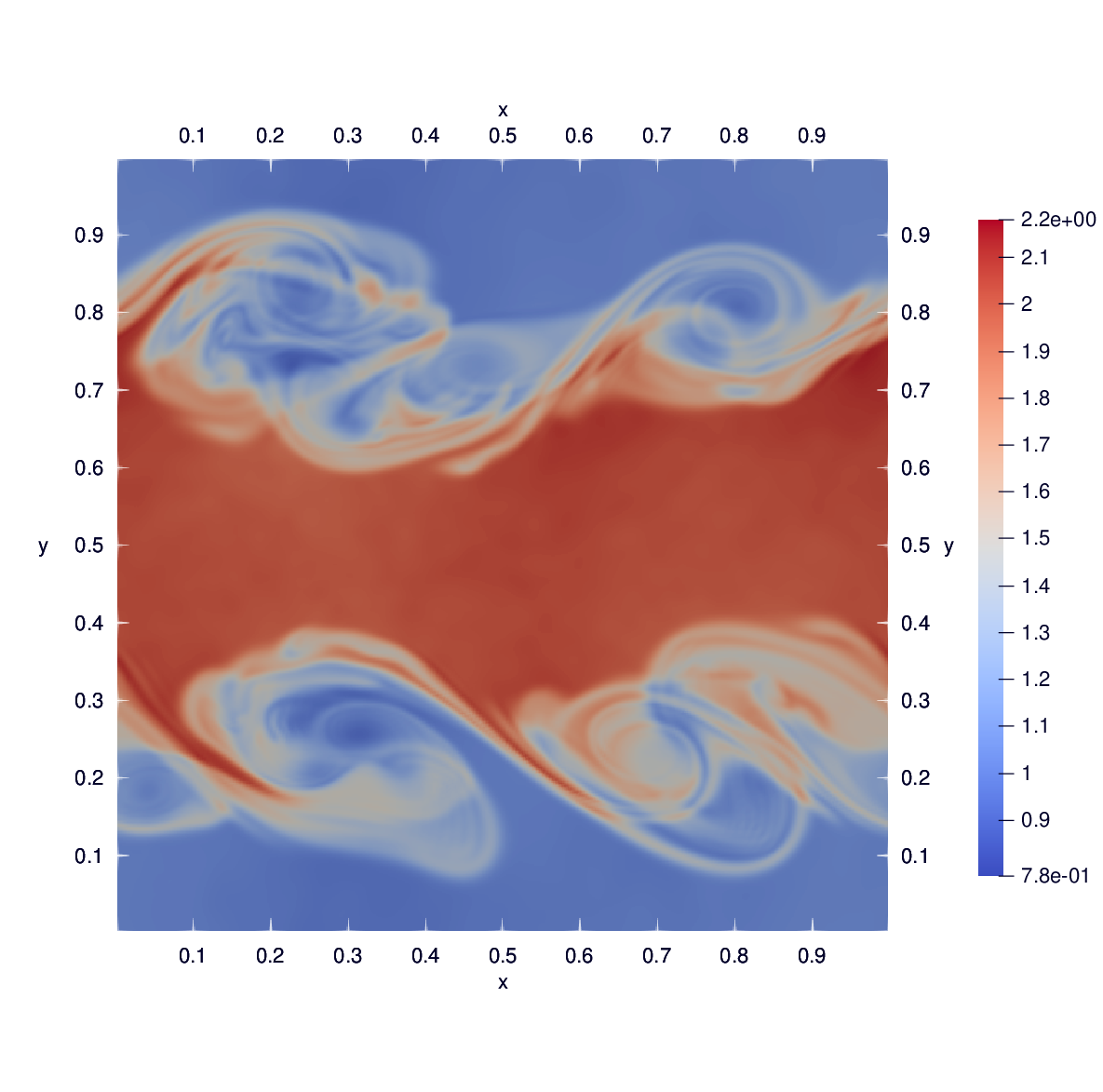}
        \vspace{-.8cm}
        \caption{Mean over 2 meshes}
        \label{fig:Cesaro_2_256}
    \end{subfigure}
    \begin{subfigure}{0.5\textwidth}
        \includegraphics[scale=0.27]{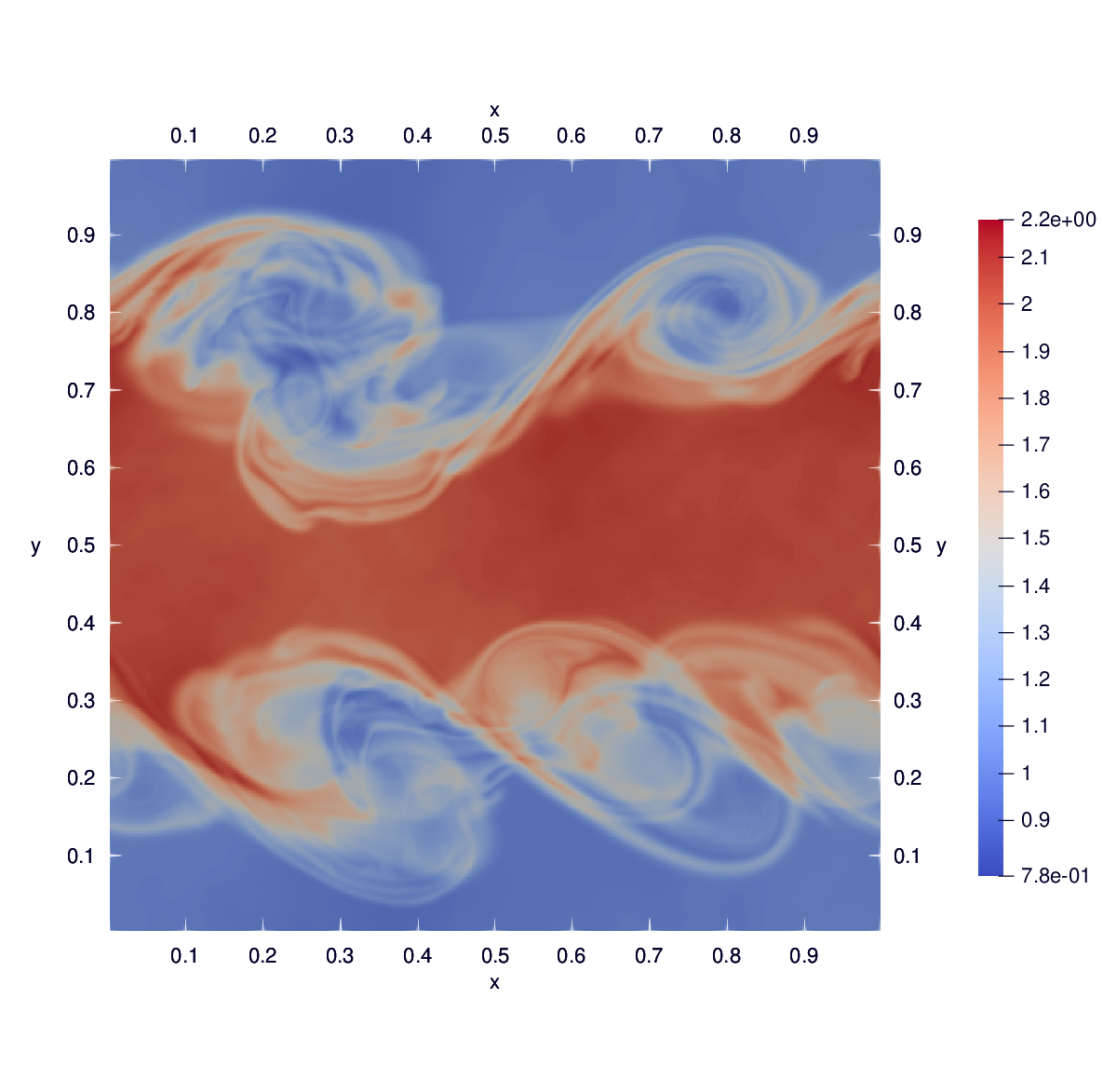}
        \vspace{-.8cm}
        \caption{Mean over 3 meshes}
        \label{fig:Cesaro_3_256}
    \end{subfigure}

    \vspace{-0.6cm}
    \begin{subfigure}{0.5\textwidth}
        \includegraphics[scale=0.27]{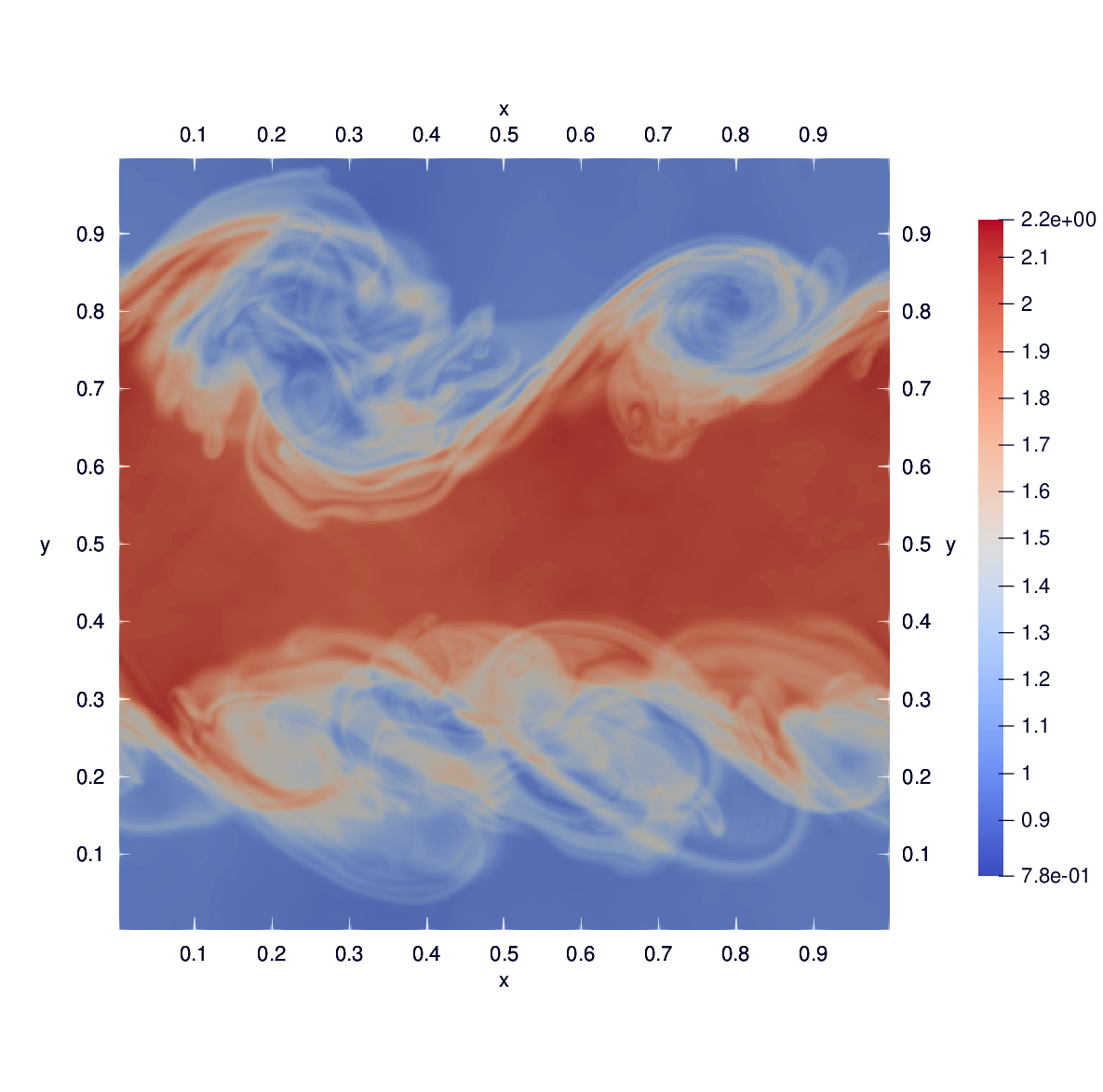}
    \vspace{-.8cm}
        \caption{Mean over 4 meshes}
        \label{fig:Cesaro_4-256}
    \end{subfigure}
    \begin{subfigure}{0.5\textwidth}
        \includegraphics[scale=0.27]{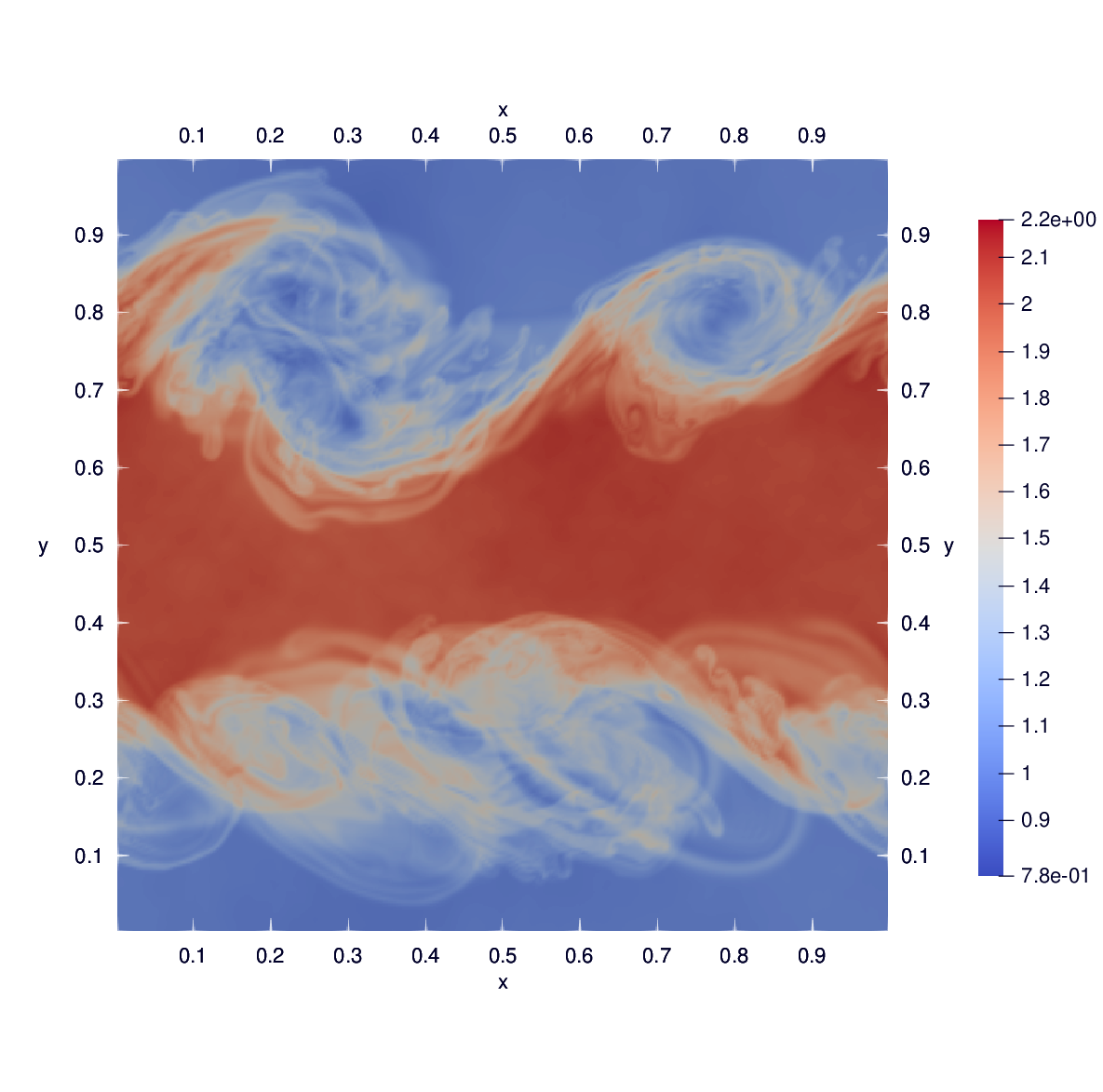}
        \vspace{-.8cm}
        \caption{Mean over 5 meshes}
        \label{fig:Cesaro_5_256}
    \end{subfigure}

    \vspace{-0.6cm}
    \centering
    \begin{subfigure}{0.5\textwidth}
        \includegraphics[scale=0.27]{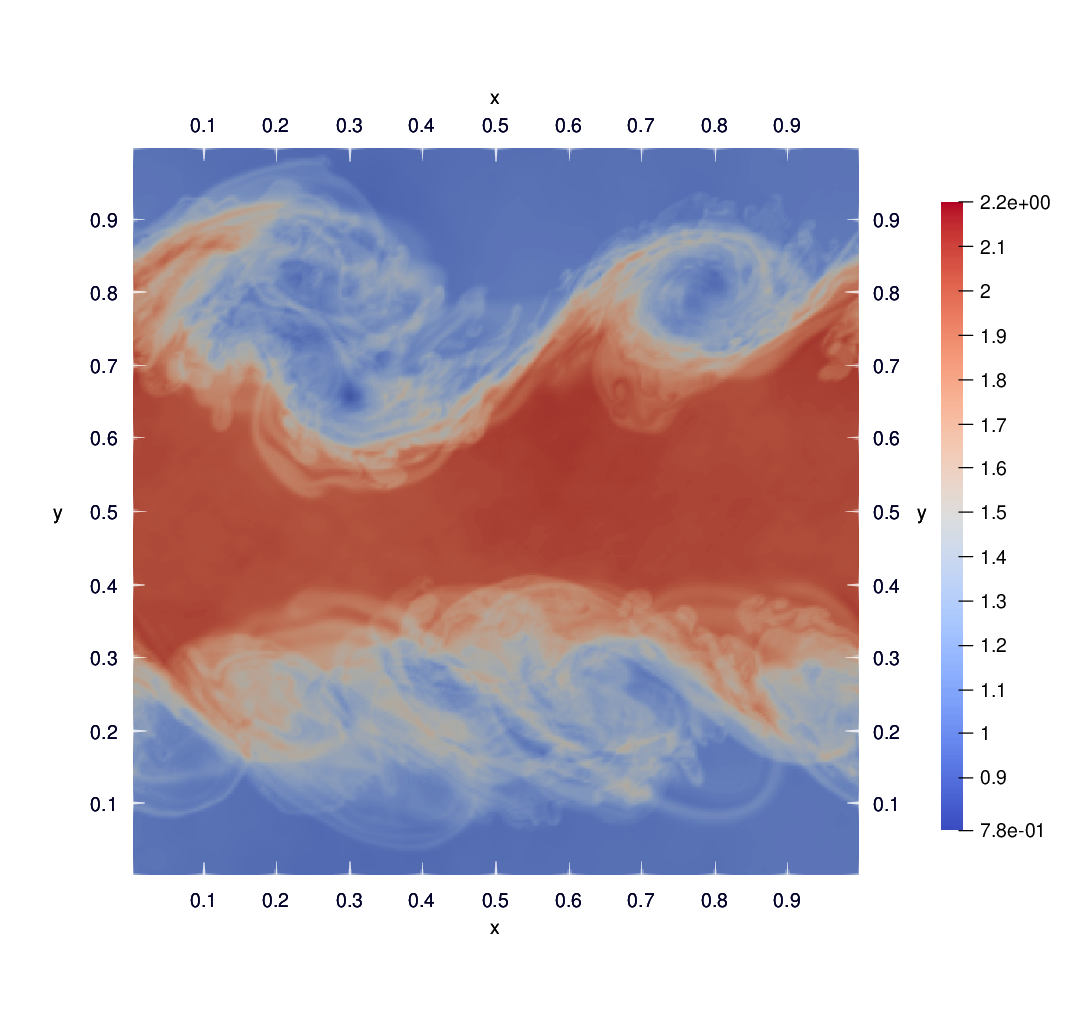}
        \vspace{-.8cm}
        \caption{{Mean over 6 meshes}}
        \label{fig:Cesaro_6_256}
    \end{subfigure}
    \caption{Kelvin-Helmholtz problem: Ces\'aro averages $\tilde{U}_{h_n}$ of the density computed with the HTC scheme at time $T=2$  on meshes with $n \times n$ cells, starting from  $n=256$.}
    \label{fig:Cesaro}
\end{figure}
\begin{figure}[!h]

    \vspace{-5mm}

    \begin{subfigure}{0.5\textwidth}
        \includegraphics[scale=0.27]{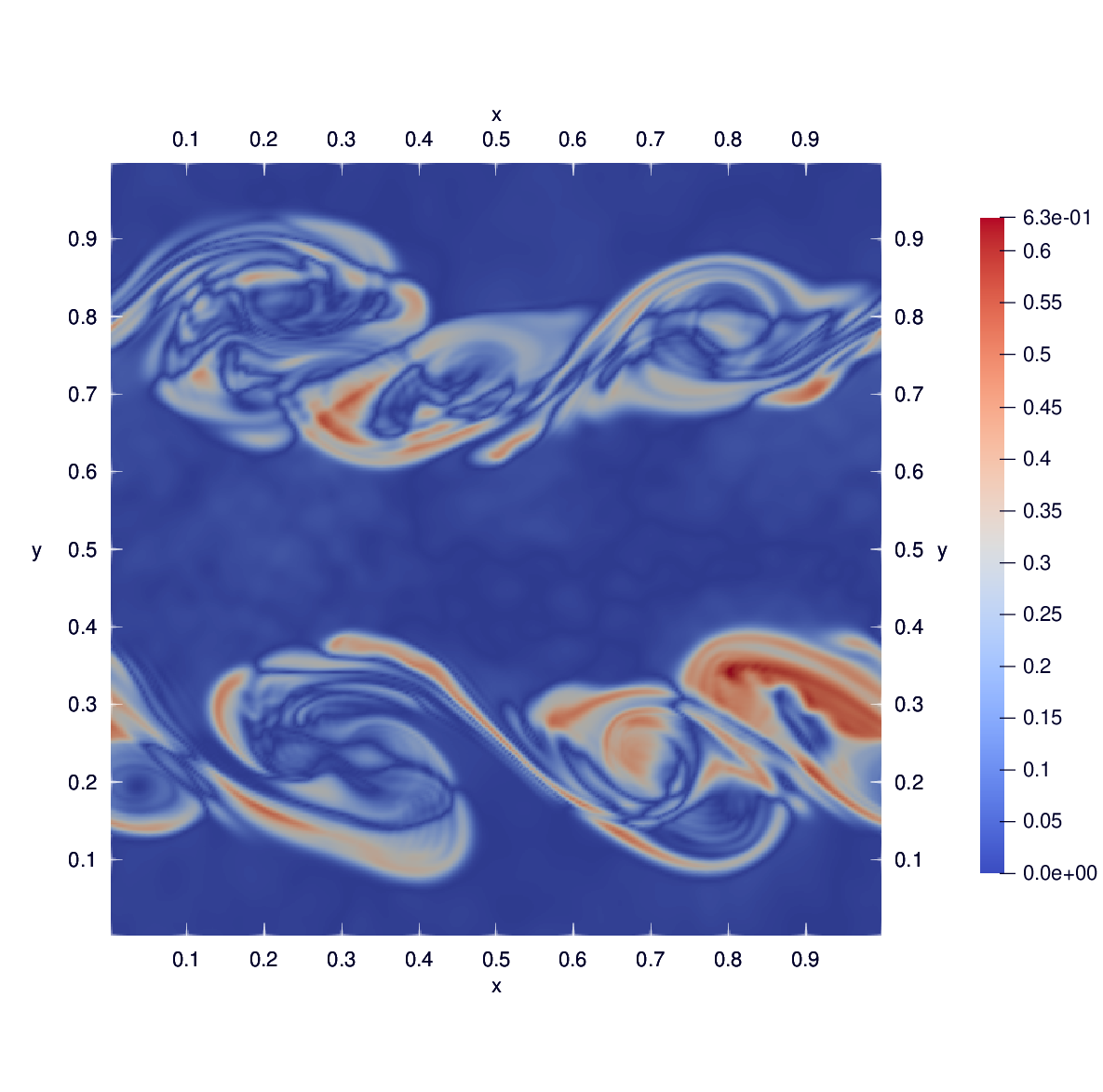}
        \vspace{-.8cm}
        \caption{First variance over 2 meshes}
        \label{fig:deviation_2_256}
    \end{subfigure}
    \begin{subfigure}{0.5\textwidth}
        \includegraphics[scale=0.27]{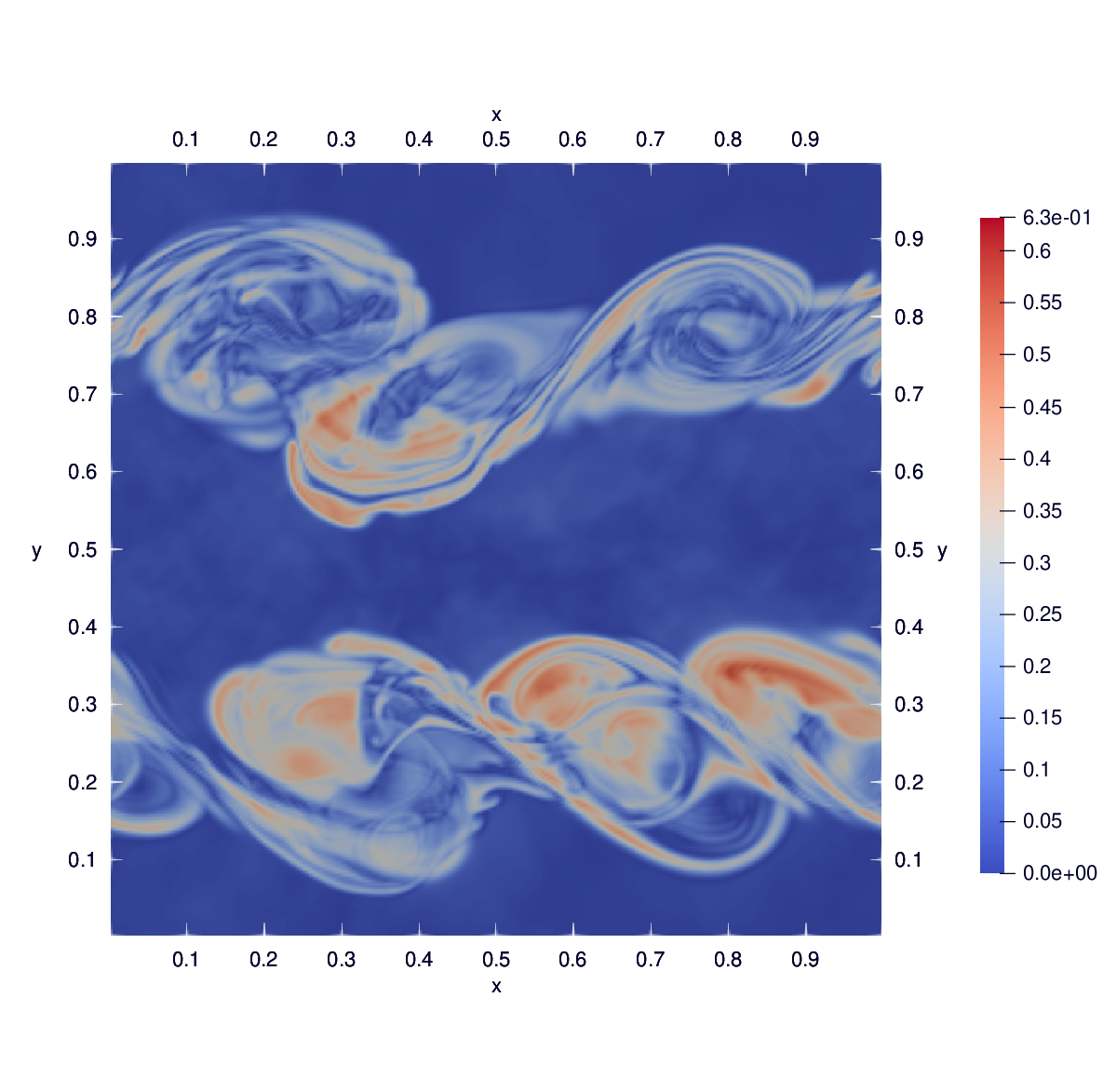}
        \vspace{-.8cm}
        \caption{First variance over 3 meshes}
        \label{fig:deviation_3_256}
    \end{subfigure}

    \vspace{-0.6cm}

    \begin{subfigure}{0.5\textwidth}
        \includegraphics[scale=0.27]{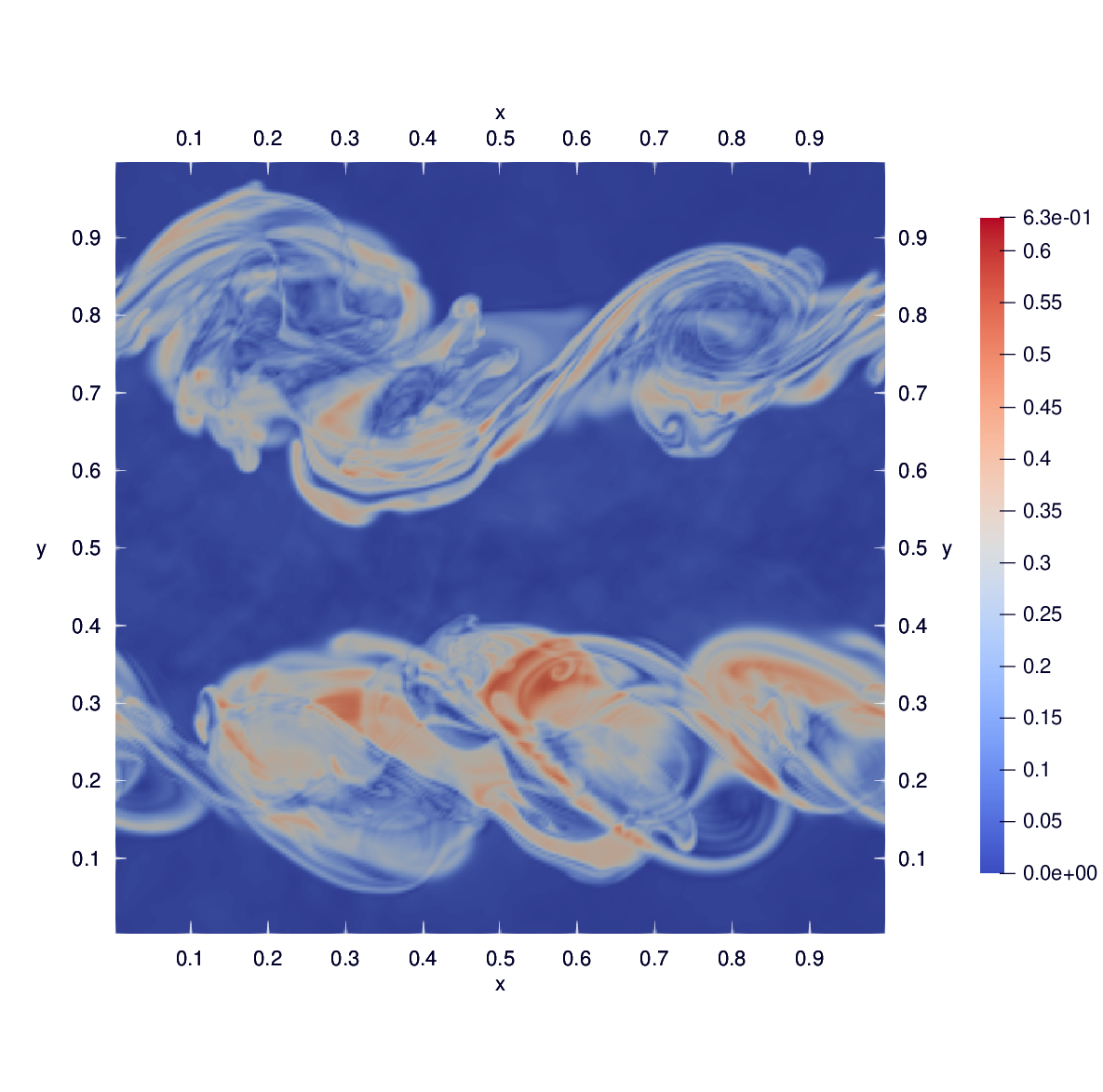}
        \vspace{-.8cm}
        \caption{First variance over 4 meshes}
        \label{fig:deviation_4_256}
    \end{subfigure}
    \begin{subfigure}{0.5\textwidth}
        \includegraphics[scale=0.27]{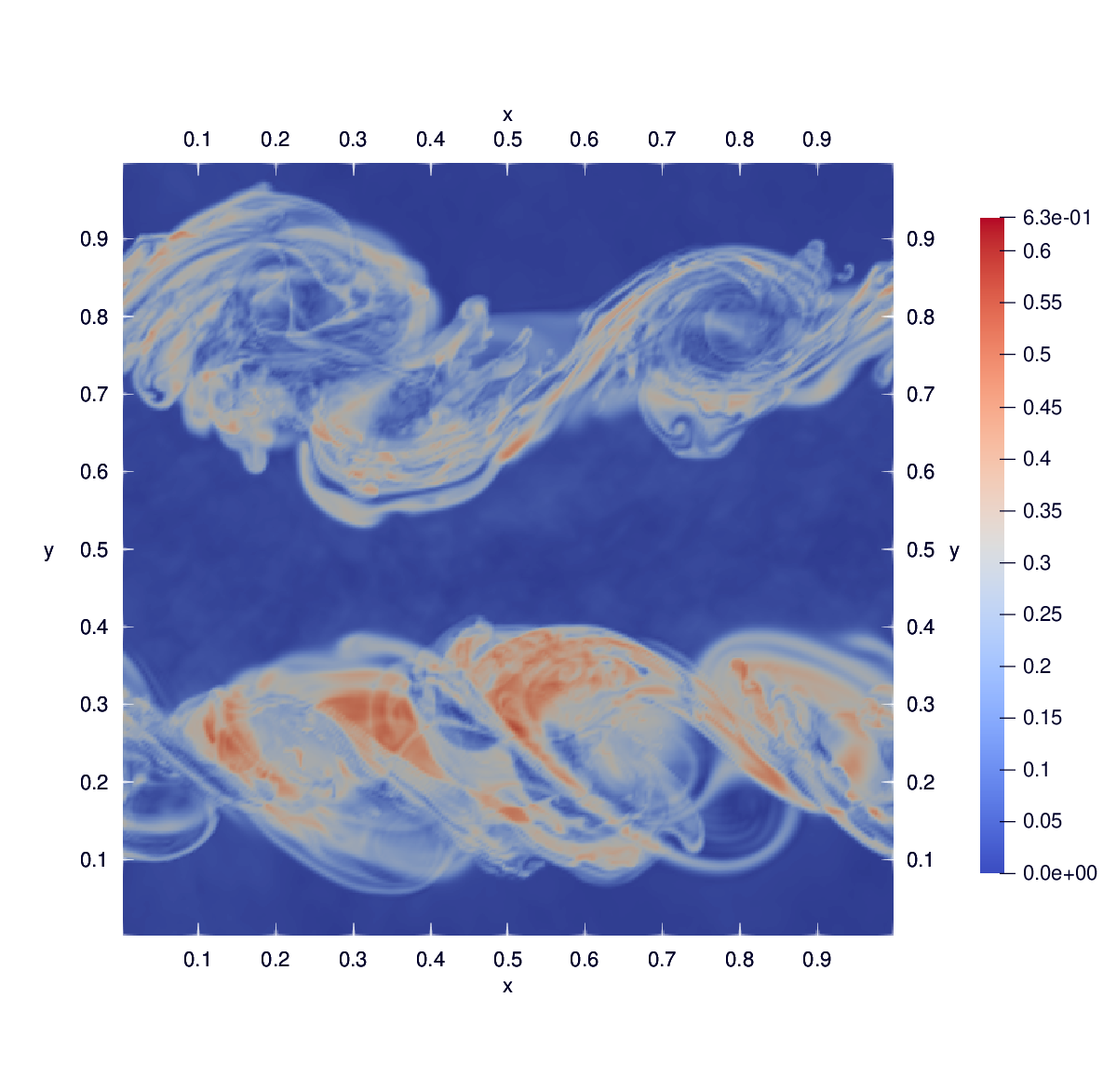}
        \vspace{-.8cm}
        \caption{First variance over 5 meshes}
        \label{fig:deviation_5_256}
    \end{subfigure}

    \vspace{-0.6cm}
    \centering
    \begin{subfigure}{0.5\textwidth}
        \includegraphics[scale=0.3,keepaspectratio]{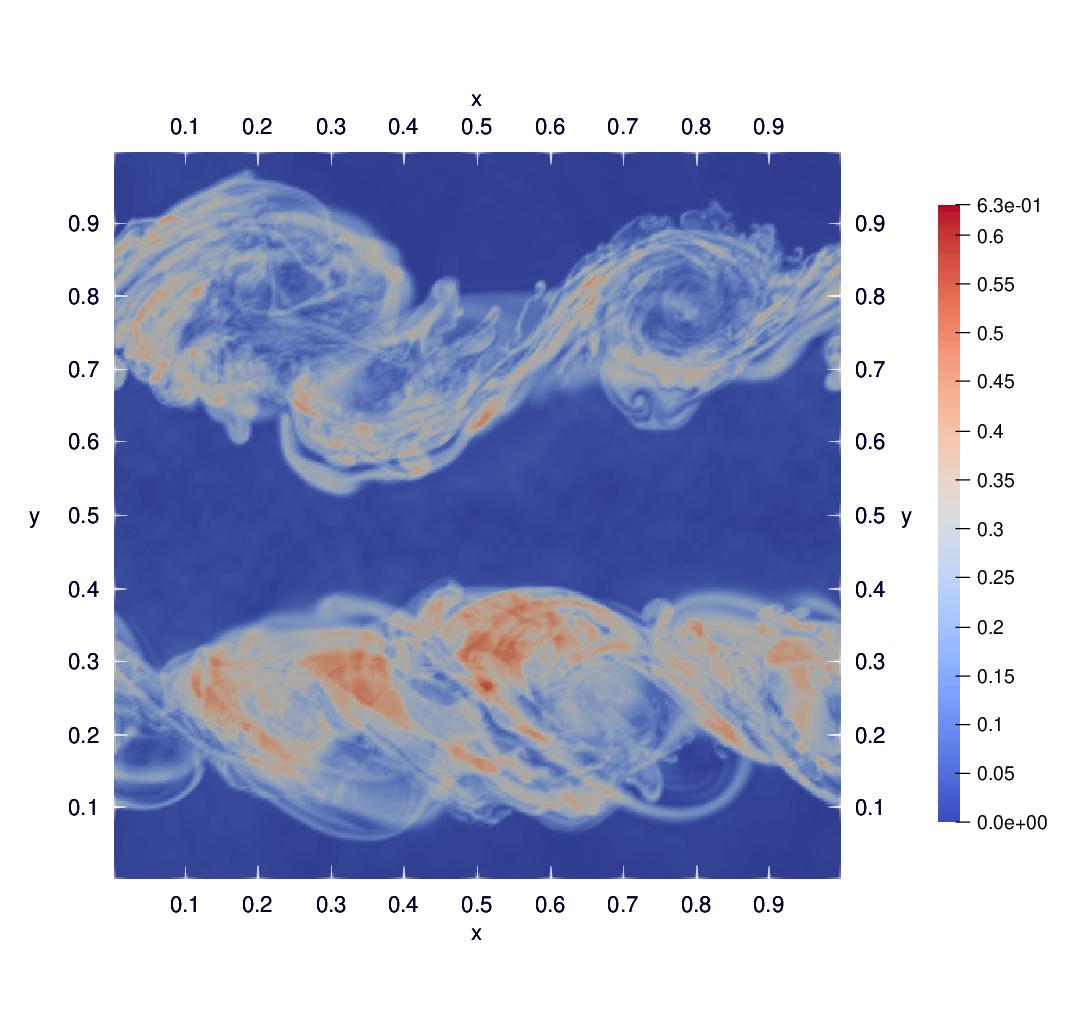}
        \vspace{-.75cm}
        \caption{{First variance over 6 meshes}}
        \label{fig:deviation_6_256}
    \end{subfigure}

    \vspace{-2mm}

    \caption{Kelvin-Helmholtz problem: First variance $U_{h_n}^\dagger$ of the density computed with the HTC scheme at time $T=2$ starting from the $n=256$ mesh.}
    \label{fig:deviation}
\end{figure}

In Table \ref{tab:EOC}, the errors and EOC for density, momentum in $x$- and $y$- direction, the total entropy and total energy are given.
As expected, the HTC scheme does not converge strongly to the reference solution computed on a grid with  $N \times N$ cells,  $N=2048$.
However, we observe the convergence of the Ces\'aro averages and the first variance.
In Figure \ref{fig:KH}, the numerical solutions for the density on  $n \times n$ meshes with $n= 256, \dots, {8192}$ cells are plotted.
We see that as the mesh is refined, smaller vortices are developed, but the structure of the solution differs and no convergence can be observed.
This is reflected also in the $E_1$-errors and corresponding EOC.
However, when looking at the solution in terms of the Ces\'aro averages given in Figure \ref{fig:Cesaro}, a clear structure can be identified
for the density when more and more meshes are taken into account.
The same behaviour can be observed in Figure \ref{fig:deviation}, where the first variance $U_{h_n}^\dagger$ is plotted.
For the remaining variables, similar results are observed and were thus omitted.
This clearly illustrates theoretical findings on the strong and weak convergence of the novel HTC scheme applied to the Euler equations.

\section{Conclusions}
\label{sec:conclusions}

We have presented the first convergence analysis of the novel HTC scheme applied to the Euler equations written as mass and momentum conservation and an entropy inequality.
A key point is the computation of time evolution of the entropy.
The structure of the system of partial differential equations allows the derivation of a positive entropy production term which is compatible with the conservation of total energy.

To achieve the convergence results we required a uniform lower bound on the density and an upper bound of the total energy, leading to crucial bounds on the temperature which was necessary to obtain boundedness of the entropy source term.
The latter implies additional stability estimates on derivatives and weak BV-estimates.
In Section \ref{sec:stability} we showed the stability of the numerical scheme which was followed by the consistency proof in Section \ref{sec:Consistency} summarized in Theorem \ref{theo: consistency}.
Finally, the convergence via dissipative weak solutions was established in Section~\ref{sec:convergence}, where the main result was presented in
Theorem~\ref{theo:DMVsolutions}.
We have showed that the approximate solutions generated by the HTC scheme converge weakly to a dissipative weak solution.
Employing dissipative weak-strong uniqueness result from \cite{Brezina2018}, the strong convergence of the HTC scheme to a strong solution of the Euler equations was established on the lifespan of the strong solution in Theorem~\ref{theo:strongconv}.
To represent the convergence to a dissipative weak solution, the notion of $\mathcal{K}$-convergence was applied.
This leads to a strong convergence of the Ces\'aro averages and the first variance in the $L^1$-norm.
Kelvin-Helmholtz instability problem with discontinuous initial data was considered for which individual realisations do not converge strongly.
However, the convergence of the Ces\'aro averages and the first variance was confirmed.
This was also reflected in the structures exhibited in  numerical solutions viewed in terms of the Ces\'aro averages and the first variance plots.

\section*{Acknowledgements}

The work of  M.L.-M. was supported by the Deutsche Forschungsgemeinschaft (DFG, German Research Foundation) - project number 233630050 - TRR 146 and
project number 525800857 -  SPP 2410 ``Hyperbolic Balance Laws: Complexity, Scales and Randomness".
She is also grateful to the  Gutenberg Research College and Mainz Institute of Multiscale Modelling  for supporting her research.
M.D. was funded by the Italian Ministry of Education, University and Research (MIUR) in the framework of the PRIN 2022 project  \textit{High order structure-preserving semi-implicit schemes for hyperbolic equations} and via the  Departments of Excellence  Initiative 2018--2027 attributed to DICAM of the University of Trento (grant L. 232/2016). M.D. was also funded by the European Union Next Generation EU projects PNRR Spoke 7 CN HPC and PNRR Spoke 7 RESTART, as well as by the European Research Council (ERC) under the European Union Horizon 2020 research and innovation programme, Grant agreement No. ERC-ADG-2021-101052956-BEYOND. Views and opinions expressed are however those of the author(s) only and do not necessarily reflect those of the European Union or the European Research Council. Neither the European Union nor the granting authority can be held responsible for them.
M.D. is member of the Gruppo Nazionale Calcolo Scientifico - Istituto Nazionale di Alta Matematica (GNCS-INdAM).

%

    \bibliographystyle{plain}
    \bibliography{lit_HTC_conv.bib}
    %

\end{document}